\documentclass[a4paper, notitlepage, 12pt]{article}

\usepackage{ucs}                
\usepackage[utf8x]{inputenc} 	

\usepackage[T1]{fontenc}     
\usepackage{lmodern}         

\usepackage[singlespacing]{setspace}

\usepackage[total={160mm,250mm}, top=25mm, left=17mm, includefoot]{geometry}

\usepackage{color}
\usepackage{lscape}
\usepackage{amsmath}
\usepackage{amssymb}
\usepackage{amsbsy}
\usepackage{amsthm}
\usepackage{wasysym}
\usepackage{float}
\usepackage{algpseudocode}
\usepackage[section]{algorithm}
\usepackage{tensor}
\usepackage[breaklinks]{hyperref}
\usepackage{authblk}

\usepackage{multirow} 
\usepackage{array}
\usepackage{graphicx}
\usepackage{caption}

\usepackage{rcs} \RCS $Revision: 2.2 $ \RCS $Date: 2013/12/17 22:30:07 $ \RCS $Author: hgm $ \RCS $RCSfile: NonLinBU.tex,v $ \RCS $Id: NonLinBU.tex,v 2.2 2013/12/17 22:30:07 hgm Exp $

\usepackage[british]{babel}

\newcommand{\citep}[1]{\cite{#1}}

\usepackage{refdef}
\usepackage{nfontdef}

\newtheorem{thm}{Theorem}

\newtheorem{prop}[thm]{Proposition}

\newcommand{\ignore}[1]{}

\newcommand{\vrho}{\varrho}
\newcommand{\vepsilon}{\varepsilon}
\newcommand{\vsigma}{\varsigma}

\newcommand{\vphi}{\varphi}

\newcommand{\vTheta}{\varTheta}

\newcommand{\vXi}{\varXi}

\newcommand{\vOmega}{\varOmega}

\newcommand{\vek}[1]{\mathchoice{\displaystyle\boldsymbol{#1}}
{\textstyle\boldsymbol{#1}}{\scriptstyle\boldsymbol{#1}}
{\scriptscriptstyle\boldsymbol{#1}}}
\newcommand{\mat}[1]{\mathchoice{\displaystyle\mathbf{#1}}
{\textstyle\mathbf{#1}}{\scriptstyle\mathbf{#1}}
{\scriptscriptstyle\mathbf{#1}}}
\newcommand{\opb}[1]{\vek{{\mathsf{#1}}}}

\newcommand{\vhat}[1]{\vek{\hat{#1}}}

\newcommand{\EXP}[1]{\mathbb{E}\left(#1\right)}

\newcommand{\divg}{\mathop{\mathrm{div}}\nolimits}

\newcommand{\dd}{\partial}
\newcommand{\di}{\mathrm{d}}

\newcommand{\ip}[2]{\langle #1, #2 \rangle}
\newcommand{\bkt}[2]{\langle #1 | #2 \rangle}
\newcommand{\ipj}[1]{\langle #1 \rangle}

\newcommand{\ns}[1]{| #1 |}
\newcommand{\nd}[1]{\| #1 \|}

\newcommand{\Hf}[2]{\tensor[^#1]{#2}{}}

\newcommand{\cov}{\mbox{cov}}

\definecolor{myred}{rgb}{1, 0.2, 0.2}

\newcommand{\authorhgm}{Hermann G. Matthies}
\newcommand{\authoral}{Alexander Litvinenko}

\newcommand{\affilwire}{Institute of Scientific Computing, \authorcr
                        Technische Universit\"at Braunschweig}

\newcommand{\thetitle}{Inverse problems and {\authorcr} uncertainty quantification}

\newcommand{\theauthor}{\authoral, \authorhgm}
\newcommand{\thesubject}{(MSC2010) 62F15, 65N21, 62P30, 60H15, 60H25, 74G75, 80A23, 74C05\\ (PACS2010) 46.65.+g, 46.35.+z, 44.10.+i\\ (ACM1998) G.1.8, G.3, J.2}
\newcommand{\thekeywords}{inverse problem, identification, uncertainty quantification}
\newcommand{\textdate}{in December 2013}

\newcommand{\thebib}{./bib}

\include{config}

\begin{document}

\title{\thetitle\thanks{Partly supported by the Deutsche
          Forschungsgemeinschaft (DFG) through SFB 880.}}

\author{\authoral}
\makeatletter
\author{\authorhgm 
\thanks{Corresponding author: D-38092 Braunschweig, 
       Germany, e-mail: \texttt{wire@tu-bs.de}}
}
\makeatother

\affil{\affilwire}

\date{}


\ignore{          


\setcounter{page}{0}
\thispagestyle{empty}
\cleardoublepage

\include{titlepage}

\newpage

\thispagestyle{empty}
\vspace*{\stretch{2}}

\begin{flushleft}
\begin{tabular}{ll}
\makeatletter
This document was created \textdate{} using \LaTeXe. \\[1cm]
\makeatother
\end{tabular}

\begin{tabular}{ll}
\begin{minipage}{6cm}
Institute of Scientific Computing\\ 
Technische Universit\"at Braunschweig\\
Hans-Sommer-Stra\ss{}e 65\\
D-38106 Braunschweig, Germany\\

\texttt{url: \url{www.wire.tu-bs.de}}\\
\makeatletter
\texttt{mail: \href{mailto:wire@tu-bs.de?subject=\thetitle}{wire@tu-bs.de}}
\makeatother
\end{minipage}
&
\begin{minipage}{2.5cm}
\vspace{-0.5cm}
\includegraphics[scale=0.34]{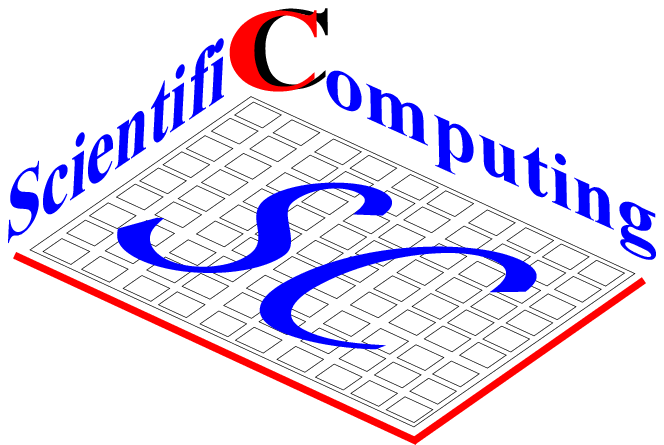}

\end{minipage}
\end{tabular}

\vspace*{\stretch{1}}

Copyright \copyright{} by \theauthor{}\\[5mm]
\end{flushleft}

This work is subject to copyright. All rights are reserved, whether the whole or part of the material is concerned, specifically the rights of translation, reprinting, reuse of illustrations, recitation, broadcasting, reproduction on microfilm or in any other way, and storage in data banks. Duplication of this publication or parts thereof is permitted in connection with reviews or scholarly analysis. Permission for use must always be obtained from the copyright holder.\\[5mm]

Alle Rechte vorbehalten, auch das des auszugsweisen Nachdrucks, der auszugsweisen oder vollständigen Wiedergabe (Photographie, Mikroskopie), der Speicherung in Datenverarbeitungsanlagen und das der Übersetzung.


}            

\maketitle

%

\begin{abstract}
In a Bayesian setting, inverse problems and uncertainty quantification (UQ)---the
propagation of uncertainty through a computational (forward)
model---are strongly connected.
In the form of conditional expectation the Bayesian update becomes computationally
attractive.  This is especially the case as together with a functional or
spectral approach for the forward UQ there is no need for time-consuming and slowly
convergent Monte Carlo sampling.  The developed sampling-free non-linear
Bayesian update is derived from the variational problem associated with conditional
expectation.  This formulation in general calls for further discretisation to
make the computation possible, and we choose a polynomial approximation.
After giving details on the actual computation in the framework of functional
or spectral approximations, we demonstrate the workings of the algorithm on a
number of examples of increasing complexity.  At last, we compare the linear
and quadratic Bayesian update on the small but taxing example of
the chaotic Lorenz 84 model, where we experiment with the influence of
different observation or measurement operators on the update.

\vspace{5mm}

\vspace{5mm}
{\noindent\textbf{Classification:} \thesubject}

\end{abstract}

%
%
%
%
%
%









%

\section{Introduction}  \label{S:intro}
In trying to predict the behaviour of physical systems, one is often
confronted with the fact that although one has a mathematical model
of the system which carries some confidence as to its fidelity, some
quantities which characterise the system may only be incompletely
known, or in other words they are uncertain.
See \cite{boulder:2011} for a synopsis on our approach to
such parametric problems.

We want to identify these parameters through
observations or measurement of the response of the system,
which can be approached in different ways.
In the mathematical description, the measurement / observation / output is
determined by the uncertain parameters, i.e.\ we have a mapping
from parameters to observations.  The problems is that usually
this mapping is not invertible, hence these inverse identification
problems are generally ill-posed.

One way to deal with this difficulty is to measure the difference
between observed and predicted system output and try to find parameters such
that this difference is minimised.  Frequently it may happen that
the parameters which realise the minimum are not unique.  In case
one wants a unique parameter, a choice has to be made, usually by
demanding additionally that some norm or similar functional of the parameters
is small as well, i.e.\ some regularity is enforced. This optimisation
approach hence leads to regularisation procedures \cite{Engl2000}.

Here we take the view that our lack of knowledge or uncertainty
of the actual value of
the parameters can be described in a \emph{Bayesian} way through a
probabilistic model \citep{jaynes03, Tarantola2004, Stuart2010}.  
The unknown parameter is then modelled as a random variable (RV)---also 
called the \emph{prior} model---and additional information on the
system through measurement or observation
changes the probabilistic description to the so-called \emph{posterior} model.
The second approach is thus a method to update the probabilistic description 
in such a way as to take account of the additional information, and the 
updated probabilistic description \emph{is} the parameter estimate,
including a probabilistic description of the remaining uncertainty.

It is well-known that such a Bayesian update is in fact closely related
to \emph{conditional expectation}
\citep{Bobrowski2006/087, Goldstein2007}, and this will be the basis
of the method presented.  For these and other probabilistic notions
see for example \citep{Papoulis1998/107} and the references therein.
As the Bayesian update may be numerically
very demanding, we show computational procedures
to accelerate this update through methods based on 
\emph{functional approximation} or \emph{spectral representation}
of stochastic problems \citep{matthies6}.  These approximations
are in the simplest case known as Wiener's so-called \emph{homogeneous}
or \emph{polynomial chaos} expansion 
\citep{Wiener1938}, which are polynomials in independent Gaussian RVs
---the `chaos'---and which
 can also be used numerically in a Galerkin procedure
\citep{ghanemSpanos91, matthiesKeese05cmame, matthies6}.  This approach
has been generalised to other types of RVs \citep{xiuKarniadakis02a}.
It is a computational variant of \emph{white noise analysis}, which means
analysis in terms of independent RVs, hence the term `white noise'
\citep{holdenEtAl96, Janson1997, hida}, see also
\citep{matthiesKeese05cmame, Roman_Sarkis_06}, and
\citep{GalvisSarkis:2012} for here relevant results on stochastic
regularity.  Here we describe computational extensions of this 
approach to the inverse problem of Bayesian updating, see also
\citep{opBvrAlHgm12, bvrAlOpHgm12-a, OpBrHgm12, BvrAkJsOpHgm11}.

To be more specific, let us consider the following situation:
we are investigating some physical system which
is modelled by an evolution equation for its state: 
\begin{equation} \label{eq:I}
\frac{\dd}{\dd t}u(t) + A(q;u(t)) = f(q;t),
\end{equation}
where $u(t) \in \C{U}$ describes the state of the system
at time $t \in [0,T]$ lying in a Hilbert space $\C{U}$ (for the sake of
simplicity), $A$ is a---possibly non-linear---operator modelling
the physics of the system, and $f\in\C{U}^*$ is some external
influence (action / excitation / loading).  The model depends on some parameters
$q \in \C{Q}$ which are uncertain and which we would thus like to identify.
To have a concrete example of \refeq{eq:I}, consider the diffusion equation
\begin{equation} \label{eq:I-c}
\frac{\dd}{\dd t}u(x,t) - \divg (\kappa(x) \nabla u(x,t)) = f(x,t),
\quad x\in\C{G},
\end{equation}
with appropriate boundary and initial conditions, where $\C{G} \subset \D{R}^n$
is a suitable domain.  The diffusing quantity is $u(x,t)$ (heat, concentration)
and the term $f(x,t)$ models sinks and sources.
Similar examples will be used for the numerical
experiments in \refS{bayes-lin} and \refS{bayes-non-lin}.
Here $\C{U} = H^1_E(\C{G})$, the subspace
of the Sobolev space $H^1(\C{G})$ satisfying the essential boundary 
conditions, and we assume that the diffusion coefficient $\kappa(x)$ is uncertain.
The parameters could be the positive diffusion coefficient
field $\kappa(x)$, but for reasons to be explained fully later we prefer
to take $q(x) = \log(\kappa(x))$, and assume $q \in \C{Q} = L_2(\C{G})$.

Our main application focus are models described by partial differential
equations (PDEs) like \refeq{eq:I-c}, and discretised for 
example by finite element procedures.
The updating methods have to be well defined and stable in a continuous
setting, as otherwise one can not guarantee numerical stability with
respect to the PDE discretisation refinement, see \citep{Stuart2010}
for a discussion of related questions.  Due to this we describe the update
before any possible discretisation in the simplest Hilbert space setting.
On the other hand
no harm will result for the basic understanding if the reader wants to 
view the occurring spaces as finite dimensional Euclidean spaces.

Now assume that we observe a function of
the state $Y(u(q),q)$, and from this observation we would like to identify
the corresponding $q$.  In the concrete example \refeq{eq:I-c} this could
be the value of $u(x_j,t)$ at some points $x_j \in \C{G}$.
This is called the \emph{inverse} problem, and
as the mapping $q\mapsto Y(q)$ is usually not invertible, the inverse
problem is \emph{ill-posed}.  Embedding this problem of finding the best
$q$ in a larger class by modelling our knowledge
about it with the help of probability theory, then in a Bayesian manner the
task becomes to estimate conditional expectations, e.g.\ see 
\cite{jaynes03, Tarantola2004, Stuart2010}
and the references therein.  The problem now is
\emph{well-posed}, but at the price of `only' obtaining probability
distributions on the possible values of $q$, which now is modelled
as a $\C{Q}$-valued random variable (RV).  On the other hand one
naturally also obtains information about the remaining uncertainty.
Predicting what the measurement
$Y(q)$ should be from some assumed $q$ is computing the \emph{forward}
problem.  The \emph{inverse} problem is then approached by comparing the
forecast from the forward problem with the actual information.

Since the parameters of the model to be estimated are uncertain, all relevant
information may be obtained via their stochastic description.
In order to extract information from the posterior, most estimates take
the form of expectations w.r.t.\ the posterior.
These expectations---mathematically integrals, numerically to be evaluated
by some quadrature rule---may be computed via asymptotic,
deterministic, or sampling methods.
  In our review of current work we
follow our recent publications \citep{opBvrAlHgm12, bvrAlOpHgm12-a, OpBrHgm12,
BvrAkJsOpHgm11}.

One often used technique is a 
Markov chain Monte Carlo (MCMC) method \cite{Madras-Fields:2002, Gamerman06},
constructed such
that the asymptotic distribution of the Markov chain is the Bayesian
posterior distribution; for further information see \citep{BvrAkJsOpHgm11}
and the references therein.

These approaches require a large number of samples in order
to obtain satisfactory results.  Here
the main idea here is to perform the Bayesian update directly on the
polynomial chaos expansion (PCE) without any sampling 
\citep{opBvrAlHgm12, bvrAlOpHgm12-a, boulder:2011, OpBrHgm12, BvrAkJsOpHgm11}.
This idea has appeared independently in \citep{Blanchard2010a}
in a simpler context, whereas
in \citep{saadGhn:2009} it appears as a variant of the Kalman filter
(e.g.\ \citep{Kalman}).
A PCE for a push-forward of the posterior
measure is constructed in \citep{moselhyYMarz:2011}.

From this short overview it becomes apparent that
the update may be seen abstractly
in two different ways.
Regarding the uncertain parameters 
\begin{equation}  \label{eq:RVq}
q: \vOmega \to \C{Q} \text{  as a RV on a probability space   }
  (\vOmega, \F{A}, \D{P})
\end{equation}
where the set of elementary events is $\vOmega$, $\F{A}$ a $\sigma$-algebra of
events, and $\D{P}$ a probability measure, one set of methods performs
the update by changing the probability measure $\D{P}$ and leaving the
mapping $q(\omega)$ as it is, whereas the other set of methods leaves the
probability measure unchanged and updates the function $q(\omega)$.
In any case, the push forward measure $q_* \D{P}$ on $\C{Q}$
defined by $q_* \D{P}(\C{R}) := \D{P}(q^{-1}(\C{R}))$ for a measurable
subset $\C{R} \subset \C{Q}$ is changed from prior to posterior.  For the
sake of simplicity we assume here that $\C{Q}$---the set containing possible
realisations of $q$---is a Hilbert space.  If the parameter $q$ is a RV,
then so is the state $u$ of the system \refeq{eq:I}.  In order to avoid a
profusion of notation, unless there is a possibility of confusion,
we will denote the random variables $q, f, u$ which now take values
in the respective spaces $\C{Q}, \C{U}^*$ and $\C{U}$ with the same symbol as
the previously deterministic quantities in \refeq{eq:I}.

In our overview on \citep{BvrAkJsOpHgm11} spectral methods in
identification problems we show that Bayesian identification
methods \citep{jaynes03, Tarantola2004, Goldstein2007, Stuart2010} are
a good way to tackle the identification problem, especially
when these latest developments in functional approximation methods
are used.  In the series of papers \citep{opBvrAlHgm12, bvrAlOpHgm12-a,
boulder:2011, OpBrHgm12, BvrAkJsOpHgm11}, Bayesian updating has been
used in a linearised form, strongly related to the Gauss-Markov
theorem \citep{Luenberger1969}, in ways very similar to the well-known
Kalman filter \citep{Kalman}.  This turns out to be a linearised version
of \emph{conditional expectation}.  Here we want to extend this
to a non-linear form, and show some examples of linear (LBU) and
non-linear (NLBU) Bayesian updates.

The organisation of the remainder of the paper is as follows: 
in \refS{bayes} we review the Bayesian update---classically defined via
conditional probabilities---and recall the link between conditional 
probability measures and conditional expectation.  We show how to
approximate this up to any desried polynomial degree, not only the
linearised version \citep{Luenberger1969, Kalman} which was used in
\citep{opBvrAlHgm12, bvrAlOpHgm12-a, boulder:2011, OpBrHgm12, BvrAkJsOpHgm11}.

The numerical realisation in terms of a functional or spectral 
approximation---here we use the well known Wiener-Hermite chaos---is
shortly sketched in \refS{num-real}.
In \refS{bayes-lin} we then show some computational
examples with the \emph{linear version (LBU)}, whereas in \refS{bayes-non-lin}
we show how to compute with the non-linear version.  Some concluding
remarks are offered in \refS{concl}.

%
%
%
%
%
%
%
%
%


%

\section{Bayesian Updating} \label{S:bayes}
In the setting of \refeq{eq:I} let us pose the following problem:
the parameters $q \in \C{Q}$ are uncertain or unknown.
By making observations $z_k$ at times $0 < t_1 < \dots < t_k \dots \in [0,T]$
one would like to infer what they are.  But we can not observe
the entity $q$ directly---like in Plato's cave allegory we can only see
a `shadow' of it, formally given by a `measurement operator'
\begin{equation}  \label{eq:iI}
Y: \C{Q} \times \C{U} \ni (q,u(t_k)) \mapsto y_k = Y(q; u(t_k)) \in \C{Y};
\end{equation}
at least this is our model of what we are measuring.
We assume that the space of possible measurements $\C{Y}$ is a vector space,  
which frequently may be regarded as finite dimensional,
as one can only observe a finite number of quantities.

Usually the observation of the `truth' $\hat{y}_k$
will deviate from what we expect to observe even if we knew the right $q$
as \refeq{eq:I} is only a \emph{model}---so there is some model error
$\epsilon$, and the measurement will be polluted by some measurement error
$\vepsilon$.
Hence we observe $z_k = \hat{y}_k + \epsilon + \vepsilon$.
From this one would like to know what $q$ and $u(t_k)$ are.  For the
sake of simplicity we will only consider one error term
$z_k = \hat{y}_k + \vepsilon$ which subsumes all the errors.

The mapping in \refeq{eq:iI} is usually not invertible and hence the problem
is called ill-posed.  One way to address this is via regularisation
(see e.g.\ \cite{Engl2000}),
but here we follow a different track.  Modelling our lack-of-knowledge
about $q$ and $u(t_k)$ in a Bayesian way \cite{Tarantola2004} by replacing them
with a $\C{Q}$- resp.\ $\C{U}$-valued random variable (RV), the problem
becomes well-posed \cite{Stuart2010}.  But of course one is looking now at the
problem of finding a probability distribution that best fits the data;
and one also obtains a probability distribution, not just \emph{one} pair
$q$ and $u(t_k)$.  Here we focus on the use of a linear Bayesian approach 
\cite{Goldstein2007} in the framework of `white noise' analysis.

We also assume that the error $\vepsilon(\omega)$ is a $\C{Y}$-valued RV.
Please observe that although $\hat{y}_k$ may be a deterministic 
quantity---the unknown `truth'---the model for the observed
quantity $z_k(\omega)= \hat{y}_k +\vepsilon_k(\omega)$ therefore 
becomes a RV as well.

The mathematical setup then is as follows: we assume that $\vOmega$
is a measure space with $\sigma$-algebra $\F{A}$ and
with a probability measure $\D{P}$, and that
$q: \vOmega \to \C{Q}$ and $u: \vOmega \to \C{U}$ are random variables (RVs).
The corresponding \emph{expectation} will be denoted by $\bar{q} = \EXP{q}
= \int_{\vOmega} q(\omega)\; \D{P}(\di \omega)$,
giving the mean $\bar{q}$ of the random variable, also denoted by 
$\ipj{q} := \bar{q}$.  The quantity $\tilde{q} := q - \bar{q}$ is the
zero-mean or fluctuating part of the RV $q$.  The covariance between
two RVs  $q$ and $u$ is denoted by $\cov_{q u} := 
\EXP{\tilde{q} \otimes \tilde{u}}$, the
expected value of the \emph{tensor} product of the fluctuating parts.
For simplicity, we shall also require $\C{Q}$ to be a Hilbert space
where each vector is a possible realisation.  This is in order to allow
to measure the distance between different $q$'s as the norm of their
difference, and to allow the operations of linear algebra to be
performed.

Bayes's theorem is commonly accepted as a consistent way to incorporate
new knowledge into a probabilistic description \cite{jaynes03, Tarantola2004}.
The elementary textbook statement of the theorem is about
conditional probabilities
\begin{equation}  \label{eq:iII}
 \D{P}(I_q|M_z) = \frac{\D{P}(M_z|I_q)}{\D{P}(M_z)}\D{P}(I_q),
\end{equation}
where $I_q$ is some subset of possible $q$'s, and $M_z$ is the information
provided by the measurement.  This becomes problematic when
the set $M_z$ has vanishing probability measure, but if all measures
involved have probability density functions (pdf), it may be formulated
as (\cite{Tarantola2004} Ch.\ 1.5)
\begin{equation}  \label{eq:iIIa}
 \pi_q(q|z) = \frac{p(z|q)}{Z_s} p_q(q),
\end{equation}
where $p_q$ is the pdf of $q$, $p(z|q)$ is the likelihood of $z=\hat{y}+\vepsilon$
given $q$, as a function of $q$ sometimes denoted by $L(q)$, and $Z_s$
(from German \emph{Zustandssumme})
is a normalising factor such that the conditional density $\pi_q(\cdot|z)$
integrates to unity.  These terms are in direct correspondence with
those in \refeq{eq:iII}.  Most computational approaches determine
the pdfs \cite{Marzouk2007, Stuart2010, Kucherova10}.  Please observe that
the model for the RV representing the error $\vepsilon(\omega)$ determines
the likelihood functions $\D{P}(M_z|I_q)$ resp.\ $p(z|q) = L(q)$.

However,  to avoid the critical cases alluded to above, Kolmogorov already
defined conditional probabilities via conditional expectation, e.g.\ see 
\cite{Bobrowski2006/087}.  Given the conditional expectation 
$\EXP{\cdot|M_z}$, the conditional
probability is easily recovered as $\D{P}(I_q|M_z) = \EXP{\chi_{I_q}|M_z}$,
where $\chi_{I_q}$ is the characteristic function of the subset $I_q$.
It may be shown that this extends the simpler formulation described by
\refeq{eq:iII} or \refeq{eq:iIIa} and is the more fundamental notion, which
we examine next.

\subsection{Conditional expectation} \label{SS:cond-expect}
The easiest point of departure for conditional expectation in our setting
is to define it not just for one piece of measurement $M_z$---which may
not even possible unambigously---but for 
sub-$\sigma$-algebras $\F{S} \subset \F{A}$.  A sub-$\sigma$-algebra $\F{S}$
is a mathematical description of a reduced possibility of randomness,
as it contains fewer events than the full algebra $\F{A}$.
The connection with a measurement $M_z$ is to take $\F{S}:=\sigma(z)$, 
the $\sigma$-algebra generated by the measurement  $z=Y(q)+\vepsilon$.
These are all events which are consistent with possible observations of
some value for $z$.

For RVs with finite variance---elements of $\C{S}:=L_2(\vOmega, \F{A},
\D{P})$---the space with the sub-$\sigma$-algebra 
$\C{S}_\infty :=L_2(\vOmega, \F{S},\D{P})$
is a closed subspace of the full space $\C{S}$ \citep{Bobrowski2006/087}.
It represents the
RVs which are possible candidates to represent the posterior, as they
are consistent with any possible observation or measurement.
For RVs in $\C{S}$ 
the conditional expectation $\EXP{\cdot | \F{S}}$ is defined as
the orthogonal projection onto the closed subspace
$\C{S}_\infty$, e.g.\ see \cite{Bobrowski2006/087}.
This allows a simple geometrical interpretation:  the difference
between the original RV and its projection has to be perpendicular
to the subspace (see \refeq{eq:iIV}), 
and the projection minimises the distance to the
original RV over the whole subspace (see \refeq{eq:iIII}).  
The square of this distance may
be interpreted as a difference in variance, tying conditional
expectation with variance minimisation; see for example
\citep{Papoulis1998/107} and the references therein for basic descriptions
of conditional expectation.

As we have to deal with $\C{Q}$-valued RVs, a bit more formalism is needed:
define the space $\E{Q} := \C{Q} \otimes \C{S}$ of $\C{Q}$-valued RVs
of finite variance, and set $\E{Q}_\infty := \C{Q} \otimes \C{S}_\infty$
for the $\C{Q}$-valued RVs with finite
variance on the sub-$\sigma$-algebra $\F{S}$, representing the new information.

The Bayesian update as conditional expectation is now simply formulated: 
\begin{equation} \label{eq:iIII}
  \EXP{q |\F{S}} := P_{\E{Q}_\infty}(q) := \text{arg min}_{\tilde{q}\in\E{Q}_\infty}
  \| q - \tilde{q} \|^2_{\E{Q}},
\end{equation}
where $P_{\E{Q}_\infty}$ is the orthogonal projector onto $\E{Q}_\infty$.  The norm
on the Hilbert tensor product in \refeq{eq:iIII} is as usually derived
from the inner product for  $p = r \otimes s \in \E{Q}: \ip{p}{p}_{\E{Q}} :=
\ip{r}{r}_{\C{Q}} \ip{s}{s}_{\C{S}}$,
so that $\|p\|_{\E{Q}} = \|r\|_{\C{Q}} \|s\|_{\C{S}}$.
Already in \citep{Kalman} it was noted that the conditional expectation
is the best estimate not only for the \emph{loss function} `distance squared',
as in \refeq{eq:iIII},
but for a much larger class of loss functions under certain distributional
constraints.  However for the above loss function this is valid without
any restrictions.

Requiring the derivative of the quadratic loss function in \refeq{eq:iIII} to
vanish---equivalently recalling the simple geometrical characterisation
mentioned just before about the orthogonality---one arrives at the well-known
orthogonality conditions.  For later reference, we collect this result in
\begin{prop}  \label{prop:cond-expect-orthog}
There is a unique minimiser to the problem in \refeq{eq:iIII}, denoted
by $\EXP{q |\F{S}} = P_{\E{Q}_\infty}(q) \in \E{Q}_\infty$, and it is
characterised by the \emph{orthogonality condition}
\begin{equation}  \label{eq:iIV}
  \forall \tilde{q} \in \E{Q}_\infty: \quad
  \ip{q  - \EXP{q|\F{S}}}{\tilde{q}}_{\E{Q}} = 0.
\end{equation}
\end{prop}
\begin{proof}
  Either by requiring the derivative of the \emph{loss function} $\| q - 
  \cdot \|^2_{\E{Q}}$ on the closed subspace $\E{Q}_\infty$ to vanish, or by
  remembering that the difference between $q$ and its best approximation
  from $\E{Q}_\infty$ is orthogonal to that subspace \citep{Luenberger1969}, one
  arrives immediately at \refeq{eq:iIV}.
  The existence and uniqueness of the best approximation follows from the fact that
  $\E{Q}_\infty = \C{Q}\otimes\C{S}_\infty$ is a closed subspace (as $\C{S}_\infty$ is
  a closed subspace), hence a closed convex set, and the loss function is
  continuous and strictly convex.  Equivalently, this says that the projection
  $P_{\E{Q}_\infty}$ is continuous and orthogonal, i.e.\ its norm is equal to unity.
  
  Alternatively, we may invoke the \emph{Lax-Milgram} lemma for \refeq{eq:iIV},
  coerciveness and continuity are trivially satisfied on the subspace
  $\E{Q}_\infty$, which is closed and hence a \emph{Hilbert} space.
\end{proof}
Let us remark that \emph{Pythagoras's} theorem implies that
\[\|P_{\E{Q}_\infty}(q)\|_{\E{Q}}^2 = \|q\|_{\E{Q}}^2 - 
             \|q - P_{\E{Q}_\infty}(q)\|_{\E{Q}}^2 .\]

To continue, note that the \emph{Doob-Dynkin} lemma
\cite{Bobrowski2006/087} assures us that if a RV like $\EXP{q|\F{S}}$ is
in the subspace $\E{Q}_\infty$, then $\EXP{q|\F{S}} = \vphi(z)$
for some $\vphi\in L_0(\E{Y};\E{Q})$, the space of measurable functions
from $\E{Y}:=\C{Y}\otimes\C{S}$ to $\E{Q}$.
We state this key fact and the resulting new 
characterisation of the conditional expectation in
\begin{prop}  \label{prop:Doob-Dynkin}
The subspace $\E{Q}_\infty = \C{Q} \otimes \C{S}_\infty$ is given by
\begin{equation}  \label{eq:iIVxx}
\E{Q}_\infty =  \overline{\textup{span}} \{\vphi \; | \; 
     \vphi(\phi,q) := \phi(Y(q)+\vepsilon);\; \phi \in L_0(\E{Y};\E{Q})
      \quad s.t.\ \vphi \in \E{Q}\}.      
\end{equation}
Finding
the conditional expectation may be seen as rephrasing
\refeq{eq:iIII} as:
\begin{equation} \label{eq:iV}
  \EXP{q |\sigma(Y)} := P_{\E{Q}_\infty}(q) = 
  \textup{arg min}_{\phi\in L_0(\E{Y};\E{Q})} \|q - \vphi(\phi,q)\|^2_{\E{Q}}.
\end{equation}
\end{prop}
\begin{proof}
 Follows directly from the Doob-Dynkin lemma.
\end{proof}

Then $q_a := P_{\E{Q}_\infty}(q)$ is called the \emph{updated}, \emph{analysis},
\emph{assimilated}, or \emph{posterior} value,
incorporating the new information.  This is the Bayesian update
expressed in terms of RVs instead of measures.  It is the estimate of
the unknown parameters $q$ after the measurement has been performed.

\subsection{Approximation of the conditional expectation} \label{SS:approx-cond-expect}
Computationally we will not be able to deal with the \emph{whole} space $\E{Q}_\infty$,
so we look at the effect of approximations.  Assume that $L_0(\E{Y};\E{Q})$ in
\refeq{eq:iV} is approximated by subspaces $L_{0,n} \subset L_0(\E{Y};\E{Q})$, where
$n\in\D{N}$ is a parameter describing the level of approximation and
$L_{0,n} \subset L_{0,m}$ if $n < m$, such that the subspaces
\begin{equation}  \label{eq:iIVxxN}
\E{Q}_n =  \textup{span} \{\vphi(\phi,q) \; | \; 
      \phi \in L_{0,n} \subset L_0(\E{Y};\E{Q}) \quad s.t.\ \vphi \in \E{Q}\}
      \subset \E{Q}_\infty   
\end{equation}
are closed and their union is dense $\overline{\bigcup_{n} \E{Q}_n} = \E{Q}_\infty$,
a consistency condition.
From \emph{C\'ea's} lemma we immediately get:
\begin{prop}  \label{prop:Cea}
Define
\begin{equation} \label{eq:poly-n} 
P_{\E{Q}_n}(q) := 
  \textup{arg min}_{\phi\in L_{0,n}} \|q - \vphi(\phi,q)\|^2_{\E{Q}}.
\end{equation}
Then the sequence $q_{a,n} := P_{\E{Q}_n}(q)$ converges to
$q_a := P_{\E{Q}_\infty}(q)$:
\begin{equation}  \label{eq:iIVxxC}
\lim_{n\to\infty} \|q_a - q_{a,n}\|^2_{\E{Q}} = 0.      
\end{equation}
\end{prop}
\begin{proof}
Well-posedness is a direct consequence of Proposition~\ref{prop:cond-expect-orthog},
and the $P_{\E{Q}_n}$ are orthogonal projections
onto the subspaces $\E{Q}_n$, hence their norms are all equal to 
unity---a stability condition.
Application of C\'ea's lemma then directly yields \refeq{eq:iIVxxC}.
\end{proof}

Here we choose the subspaces of polynomials up to degree $n$ for the purpose
of approximation, i.e.\
\[ \E{Q}_n := \textup{span}\{ \vphi \in \E{Q}\; | \; \vphi(\psi_n,q),
  \; \psi_n \textup{ a } n^{\textup{th}} \textup{ degree polynomial}\}, \]
and we remark that in case $\C{Y}$ is finite-dimensional---the usual case---then
the space of $n^{\textup{th}}$ degree polynomials is a closed space.
We may write this as
\begin{equation}  \label{eq:n-deg-pol}
 \psi_n(z) := \Hf{0}{H} + \Hf{1}{H}z + \dotsb 
+ \Hf{k}{H}z^{\vee k}+\dotsb + \Hf{n}{H}z^{\vee n},
\end{equation}
where $\Hf{k}{H} \in \E{L}^k_s(\E{Y},\E{Q})$ is symmetric
and $k$-linear; and $z^{\vee k} := \overbrace{z \vee \ldots \vee z}^{k}
:= \text{Sym}(z^{\otimes k})$ is the symmetric tensor product of
the $z$'s taken $k$ times with itself.  Let us remark here that the form of
\refeq{eq:n-deg-pol}, given in monomials, is numerically not a good form---except
for very low $n$---and straightforward use in computations is not recommended.
The relation \refeq{eq:n-deg-pol} could be re-written in some orthogonal 
polynomials---or in fact any other system of multi-variate functions; this
generalisation will be published elsewhere.
For the sake of conceptual simplicity, we stay wtih \refeq{eq:n-deg-pol} and
then have that
\begin{equation}  \label{eq:n-deg-pol-q}
 q_{a,n}(\Hf{0}{H},\dots,\Hf{n}{H}) := \psi_n(z) := \Hf{0}{H} + \dotsb 
+ \dotsb + \Hf{n}{H}z^{\vee n}
\end{equation}
is a function of the maps $\Hf{k}{H}$.  The stationarity or orthogonality
condition \refeq{eq:iIV} can then be written in terms of the $\Hf{k}{H}$.
We need the following abbreviations for any $k,\ell\in\D{N}_0$:
\[\ipj{p \otimes v^{\vee k}} := 
     \EXP{p \otimes v^{\vee k}} = \int_{\vOmega} p(\omega)\otimes
  v(\omega)^{\vee k} \, \D{P}(\di \omega) \]
and
\[ \Hf{k}{H}\ipj{z^{\vee(\ell+k)}} := 
\ipj{z^{\vee\ell}\vee (\Hf{k}{H}z^{\vee k)}}=
\EXP{z^{\vee\ell}\vee (\Hf{k}{H}z^{\vee k)}} .\]
We may then characterise the $\Hf{k}{H}$ in the following way:
\begin{thm} \label{T:n-vers}
With $q_{a,n}$ given by \refeq{eq:n-deg-pol-q}, the stationarity condition
\refeq{eq:iIV} becomes for any $n\in\D{N}_0$ ($\delta_{(\Hf{\ell}{H})}$ the
G\^ateaux derivative w.r.t.\ $\Hf{\ell}{H}$):
\begin{equation}  \label{eq:stat-H}
\forall \ell=0,\dotsc,n: \quad 
 \delta_{(\Hf{\ell}{H})}\;\|q - q_a(\Hf{0}{H},\dotsc,\Hf{n}{H}))\|_{\E{Q}}^2 = 0,
\end{equation}
which determine the $\Hf{k}{H}$ and may be concisely written as
\begin{equation}  \label{eq:cond-H}
\forall \ell=0,\dotsc,n:\quad \sum_{k=0}^n \Hf{k}{H}\ipj{z^{\vee(\ell+k)}} = 
      \ipj{q\otimes z^{\vee\ell}}.
\end{equation}
The Hankel operator matrix $(\ipj{z^{\vee(\ell+k)}})_{\ell,k} $ in the
linear equations \refeq{eq:cond-H} is symmetric and positive definite, hence
the system \refeq{eq:cond-H} has a unique solution.
\end{thm}
\begin{proof}
The relation \refeq{eq:cond-H} is the result of straightforward differentiation
in \refeq{eq:stat-H} (and division by $2$), and may be written in more detail as:
\begin{alignat*}{6}
&\ell = 0:\; \Hf{0}{H} &\dotsb  &+ \Hf{k}{H}\ipj{z^{\vee k}}
   &\dotsb &+ \Hf{n}{H}\ipj{z^{\vee n}}  =& \ipj{q}, \\
&\ell = 1: \; \Hf{0}{H}\ipj{z} &\dotsb  &+ 
   \Hf{k}{H}\ipj{z^{\vee (1+k)}} &\dotsb &+ \Hf{n}{H}\ipj{z^{\vee (1+n)}}
     =& \ipj{q \otimes z}, \\
&\vdots &\dotso   & & \vdots & & \vdots \\
&\ell = n:\; \Hf{0}{H}\ipj{z^{\vee n}} &\dotsb  &+ 
    \Hf{k}{H}\ipj{z^{\vee (n+k)}} &\dotsb &+  \Hf{n}{H}\ipj{z^{\vee 2n}} 
    =& \ipj{q \otimes z^{\vee n}}.
\end{alignat*}
Symmetry of the operator matrix is obvious---the $\ipj{z^{\vee k}}$ are the
coefficients---and positive definiteness follows
easily from the fact that it is the gradient of the functional in \refeq{eq:stat-H},
which is strictly convex.
\end{proof}

A la Penrose in \emph{`symbolic index'} notation---or the reader may just think of
indices in a finite dimensional space with orthonormal basis---the 
system \refeq{eq:stat-H} can be given yet another form:
denote in symbolic index notation $q = (q^m), z = (z^\jmath)$, and $\Hf{k}{H} = 
  (\tensor*[^k]{H}{^m_{\jmath_1}_{\dotso}_{\jmath_k}})$, then \refeq{eq:cond-H}
becomes, with the use of the Einstein convention of summation (a tensor contraction)
over repeated indices, and with the symmetry explicitly indicated:
\begin{multline}  \label{eq:symbolic}
\forall \ell=0,\dotsc,n; \;  \jmath_1 \le \ldots \le \jmath_{\ell}\le \ldots \le
 \jmath_{\ell + k}\le \ldots \le \jmath_{\ell+n}:\\
\ipj{z^{\jmath_1}\dotsm z^{\jmath_\ell}} \,
   (\tensor*[^0]{H}{^m}) + \dotsb +
  \ipj{z^{\jmath_1}\dotsm z^{\jmath_{\ell+1}}\dotsm z^{\jmath_{\ell+k}}} \,
     (\tensor*[^k]{H}{^m_{\jmath_{\ell+1}}_{\dotso}_{\jmath_{\ell+k}}}) +\\ 
     \dotsb +  \ipj{z^{\jmath_1}\dotsm z^{\jmath_{\ell+1}}
     \dotsm z^{\jmath_{\ell+n}}} \,
   (\tensor*[^n]{H}{^m_{\jmath_{\ell+1}}_{\dotso}_{\jmath_{\ell+n}}}) =
  \ipj{q^m z^{\jmath_1}\dotsm z^{\jmath_\ell}}.
\end{multline}
We see in this representation that the matrix does \emph{not} depend on $m$---it is 
identically \emph{block diagonal} after appropriate reordering, which makes the solution
of \refeq{eq:cond-H} or \refeq{eq:symbolic} much easier.

Some special cases are:
for $n=0$---\emph{constant} functions, we do not use any information from the
measurement---we have from \refeq{eq:cond-H} or \refeq{eq:symbolic} 
$q_a = \Hf{0}{H} = \ipj{q} = \EXP{q}$.  One could argue
that this is the best approximation to $q$ in absence of any further information.

The case $n=1$ in \refeq{eq:cond-H} or \refeq{eq:symbolic} is more interesting,
allowing up to \emph{linear} terms:
\begin{alignat*}{5}
&\Hf{0}{H} &+ &\Hf{1}{H}\ipj{z}  &= &\ipj{q} \\
&\Hf{0}{H}\ipj{z} &+ &\Hf{1}{H}\ipj{z \otimes z} &= &\ipj{q \otimes z}.
\end{alignat*}
Remembering that $[\cov_{qz}] = \ipj{q \otimes z}-\ipj{q}\otimes\ipj{z}$ and
analogous for  $[\cov_{zz}]$, one obtains by tensor multiplication
in \refeq{eq:symb-elim-1} with $\ipj{z}$ and
symbolic Gaussian elimination the \refeq{eq:symb-elim-2}.
\begin{alignat}{3} \label{eq:symb-elim-1}
\Hf{0}{H} &= &\ipj{q} &- \Hf{1}{H}\ipj{z}\\
\Hf{1}{H}(\ipj{z \otimes z}-\ipj{z}\otimes\ipj{z})= \Hf{1}{H}[\cov_{zz}]  &= 
    &\ipj{q \otimes z} &- \ipj{q}\otimes\ipj{z} = [\cov_{qz}]. \label{eq:symb-elim-2}
\end{alignat}
This gives 
\begin{align} \label{eq:LBU-n1-1H}
\Hf{1}{H} &= [\cov_{qz}][\cov_{zz}]^{-1} =: K\\
\Hf{0}{H} &= \ipj{q} - [\cov_{qz}][\cov_{zz}]^{-1} \ipj{z}. \label{eq:LBU-n1-0H}
\end{align}
where $K$ in \refeq{eq:LBU-n1-1H} is the well-known \emph{Kalman} gain operator 
\citep{Kalman}, so that finally
\begin{equation}  \label{eq:LBU-n1}
q_a = \Hf{0}{H} + \Hf{1}{H}z = 
      \ipj{q} + [\cov_{qz}][\cov_{zz}]^{-1}(z - \ipj{z})=\ipj{q} + K(z - \ipj{z}).
\end{equation}
This is called the \emph{linear} Bayesian update (LBU).
It is important to see \refeq{eq:LBU-n1} as a symbolic expression, especially
the inverse $[\cov_{zz}]^{-1}$ indicated there should not really be computed,
especially when $[\cov_{zz}]$ is ill-conditioned or close to singular.  The
inverse can in that case be replaced by the \emph{pseudo-inverse}, or rather
the computation of $K$, which is in linear algebra terms a \emph{least-squares}
approximation, should be done with orthogonal transformations and not by
elimination.  We will not dwell on these well-known matters here. 

The case $n=2$ can still be solved symbolically, the system to be solved is
from \refeq{eq:cond-H} or \refeq{eq:symbolic}:
\begin{alignat*}{7}
&\Hf{0}{H} &+ &\Hf{1}{H}\ipj{z} &+ &\Hf{2}{H}\ipj{z^{\otimes 2}}  &= &\ipj{q} \\
&\Hf{0}{H}\ipj{z} &+ &\Hf{1}{H}\ipj{z^{\otimes 2}} &+ &\Hf{2}{H}\ipj{z^{\otimes 3}}&
= &\ipj{q \otimes z}\\
&\Hf{0}{H}\ipj{z^{\otimes 2}} &+ &\Hf{1}{H}\ipj{z^{\otimes 3}} &+ 
&\Hf{2}{H}\ipj{z^{\otimes 4}}&= &\ipj{q \otimes z^{\otimes 2}}.
\end{alignat*}
After some symbolic elimination steps one obtains
\begin{alignat*}{7}
&\Hf{0}{H} &+ &\Hf{1}{H}\ipj{z} &+ &\Hf{2}{H}\ipj{z^{\otimes 2}}  &= &\ipj{q} \\
&0 &+ &\Hf{1}{H} &+ &\Hf{2}{H}\; \mat{F} &= &K\\
&0 &+ &0 &+ &\Hf{2}{H}\; \opb{G} &= &\vek{E},
\end{alignat*}
with the Kalman gain operator $K\in(\C{Q}\otimes\C{Y})^*$ from \refeq{eq:LBU-n1-1H},
the third order tensors $\mat{F}\in(\C{Y}^{\otimes 3})^*$ given in \refeq{eq:NLBU-n2-F},
and $\vek{E}\in (\C{Q}\otimes\C{Y}^{\otimes 2})^*$ given in \refeq{eq:NLBU-n2-E},
and the fourth order tensor $\opb{G}\in(\C{Y}^{\otimes 4})^*$ given
in \refeq{eq:NLBU-n2-G}:
\begin{align} \label{eq:NLBU-n2-F}
\mat{F} &= \left(\ipj{z^{\otimes 3}} - \ipj{z^{\otimes 2}} \otimes \ipj{z}\right) \cdot
 [\cov_{zz}]^{-1}, \\  \label{eq:NLBU-n2-E}
 \vek{E} &= \ipj{q \otimes z^{\otimes 2}} - \ipj{q}\otimes\ipj{z^{\otimes 2}} -
      K \cdot \left(\ipj{z^{\otimes 3}} - 
     \ipj{z} \otimes \ipj{z^{\otimes 2}} \right)\\  \label{eq:NLBU-n2-G}
\opb{G} &= \left(\ipj{z^{\otimes 4}}-\ipj{z^{\otimes 2}}^{\otimes 2}\right) - 
     \mat{F} \cdot \left(\ipj{z^{\otimes 3}} - 
     \ipj{z} \otimes \ipj{z^{\otimes 2}} \right),
\end{align}
where the single central dot `$\cdot$' denotes as usual a contraction over
the appropriate indices, and a colon `:' a double contraction.  From this one easily obtains the solution
\begin{align} \label{eq:NLBU-n2-2H}
 \Hf{2}{H}&= \vek{E}:\opb{G}^{-1} \\  \label{eq:NLBU-n2-1H}
 \Hf{1}{H} &= K - \vek{E}:\opb{G}^{-1} : \mat{F} \\  \label{eq:NLBU-n2-0H}
 \Hf{0}{H} &= \ipj{q} - (K-\vek{E}:\opb{G}^{-1} : \mat{F}) \cdot \ipj{z} 
        -  \vek{E}:\opb{G}^{-1} : \ipj{z^{\otimes 2}}.
\end{align}

\subsection{Prior information and mappings} \label{SS:prior}
In case one has \emph{prior} information \(\E{Q}_f\) and a prior 
estimate $q_f(\omega)$ (forecast), and a new measurement $z$ comes in
generating via  \(\sigma(z)\) a subspace \(\E{Q}_y \subset \E{Q}\),
one now needs a projection onto \(\E{Q}_a=\E{Q}_f + \E{Q}_y \),
with reformulation as an orthogonal direct sum
\[\E{Q}_a= \E{Q}_f + \E{Q}_y = \E{Q}_f \oplus (\E{Q}_y \cap \E{Q}_f^\perp)
=\E{Q}_f \oplus \E{Q}_\infty,\]
in order not to update twice with the nonzero part of $\E{Q}_y \cap \E{Q}_f$.

 The \emph{update} / \emph{conditional expectation} / 
\emph{assimilated} value is 
\[q_a = q_f + P_{\E{Q}_\infty} q = q_f + q_\infty,\]
where $q_\infty$ is the \emph{innovation}, the orthogonal projection onto
$\E{Q}_\infty$.  This is reminiscent of
\refeq{eq:LBU-n1}, where the term $\ipj{q}$ may be regarded as the 
\emph{prior} information (before any measurement is performed) and replaced
here by $q_f$, and the innovation there is $K(z - \ipj{z})$, which is
here represented by $q_\infty$.

The $n=1$ version of \refT{n-vers} is well-known,
and in conjunction with what was just stated about prior information is
of considerable practical importance; it is an extention of the
\emph{Kalman filter} \citep{Kalman, Luenberger1969, Papoulis1998/107}.
We rephrase this generalisation of the well-known \emph{Gauss-Markov}
theorem from \citep{Luenberger1969}  Chapter 4.6, Theorem 3:
\begin{thm}  \label{T:n-1-vers}
The update $q_{a,1}$, minimising $\|q - \cdot\|^2_{\E{Q}}$ over all elements generated
by affine mappings (the up to $n=1$ case of \refT{n-vers}) of
the measurement in the case with prior information $q_f$ and predicted
measurement $z = Y(q_f) + \vepsilon$ is given from the observation $\hat{z}$ by
\begin{equation}   \label{eq:iVII}
  q_{a,1} = q_f + K(\hat{z} - z),
\end{equation}
where the notation $q_{a,1}$ is according to \refeq{eq:n-deg-pol-q},
and the operator $K$ is the Kalman gain from \refeq{eq:LBU-n1}. 
\end{thm}
This update is in some ways very similar to the
`Bayes linear' approach \cite{Goldstein2007}.  We point out that
$q_a$ and $q_f$ are RVs, i.e.\ this is an equation in $\E{Q} = 
\C{Q} \otimes \C{S}$, whereas the traditional Kalman filter---which
looks superficially just like \refeq{eq:iVII} is an equation in $\C{Q}$.
Observe that $z = Y(q_f) + \vepsilon$
and that the error term is a RV.  Hence the quantity $z$ is an RV, and
\refeq{eq:iVII} is an equation between RVs.  If the mean is taken in
\refeq{eq:iVII}, one obtains the familiar Kalman filter formula \cite{Kalman}
for the update of the mean, and one may show \citep{opBvrAlHgm12}
that \refeq{eq:iVII} also contains the Kalman update for the covariance,
i.e.\ the Kalman filter is a low-order part of \refeq{eq:iVII}.
The computational strategy is now to replace and approximte the---only abstractly
given---computation of $q_a$ by the practically possible  calculation
of $q_{a,n}$ as in \refeq{eq:n-deg-pol-q}.  This means that we 
approximate $q_a$ by $q_{a,n}$ by using $\E{Q}_n \subset \E{Q}_\infty$,
and rely on Proposition~\ref{prop:Cea}.  This corresponds to some loss of
information from the measurement, but yields a managable computation.
If the assumptions of Proposition~\ref{prop:Cea} are satisfied, then one
can expect for $n$ large enough that the terms in \refeq{eq:n-deg-pol-q}
converge to zero, thus providing an error indicator on when a sufficient
accuracy has been reached.

In case the space generated by the measurements is not dense in $\E{Q}$
a residual error will thus remain, as the measurements do not contain enough
information to resolve our lack of knowledge about $q$.  Anyway,
finding $q$ is limited by the presence of the error $\vepsilon$, 
as obviously the error influences the update in \refeq{eq:iVII}.
If the measurement operator is approximated in some way---as it will be
in the computational examples to follow---this will introduce a new error,
further limiting the resolution.

It is maybe worthwhile to pursue the following idea:  The mapping we try
to approximate $q \mapsto P_{\E{Q}_\infty}(q)$ is an orthogonal projection,
hence linear.  This carries with it several suggestions on how to
change the update process.

Given a couple $q$ and  measurement operator $y=Y(q)$, one may change
the arrangement by mappings of $q$ or $y$.  On one hand one may consider
a---preferrably---injective map $\vTheta:\C{T}\to\C{Q}$ and choose 
$p=\vTheta^{-1}(q)\in\C{T}$ as parameter,
the measurement operator is then $\tilde{Y}(p) = Y(\vTheta(p))$.  This may
be useful as we want to perform essentially linear operations on $q$ like
the above mentioned projection and linear approximations to it, and
if the set where $q$ `lives' is not a linear set, this is problematic.
We will come across this example in \refS{bayes-lin}, where $q$ is positive---or
a symmetric positive definite tensor---and hence `lives' on an \emph{open} cone
in a vector space.  There we will choose $\vTheta = \exp$, and in
\citep{BvrAkJsOpHgm11} we give some arguments why this may be meaningful,
as this transformation puts us in the tangent space of the positive cone,
which is a \emph{linear} space.

On the other hand looking again at a given pair $q$ and $Y(q)$, the linear map
$q \mapsto P_{\E{Q}_\infty}(q)$ is approximated by $\psi_n(y) = \psi_n(Y(q))$
from \refeq{eq:n-deg-pol}---neglecting measurement error for the moment.
This means that when $Y$ is nonlinear in
$q$, the update map $\psi_n$ from \refeq{eq:n-deg-pol} has to somehow
`straighten' the nonlinearity out.  This opens the possibility to make
the update `easier': we update not from $Y(q)$, but from $\vXi(Y(q))$,
where $\vXi: \C{Q}\to \C{X}$ is chosen so that the composition $\Xi\circ Y$
is `less nonlinear'.  This means that in the computation of $\psi_n$,
we try to minimise the error of $\psi_n \circ \vXi$.
 Finding a suitable $\vXi$---in some way an `inverse'
of $Y$---is not easy.  Anyway, some preliminary examples where $\vXi$ has been
chosen heuristically were very promising and will be reported elsewhere.

If the mapping $Y$ is not injective, then
of course this can not be `ironed out' by any mapping $\vXi$, as we would need
to undo the loss of information from $Y$ being not injective---another sign
of ill-posedness.  The mapping $\vXi$ could be speculatively made into a set-valued
mapping to achieve this, but we would have, for a certain $y$, to find all 
$q\in Y^{-1}(y)$ to construct $\vXi$ such that it distinguishes them, not an easy task.

We close this section by pointing out a little example connected to these
considerations---suggested to us by \citep{bSprungk13}---which is a bit
disturbing and shows the possible problems involved and that
one has to be a bit careful:  Assume that $q_f=\theta$ is a single centred
Gaussian variable with variance $\vsigma^2$, and 
that the measurement operator is $Y(q)=q^2$,
i.e. all information about the sign is lost.  Assume that $\vepsilon = z - Y(q)$
is independent of $q$ and also centred.  Taking first the linear Bayesian update
(LBU) from \refT{n-1-vers} defined in \refeq{eq:LBU-n1-1H},  we have that---as 
$\EXP{q}=0$ and $\EXP{q^2 } = \vsigma^2 $
\[ [\cov_{q,z}] = [\cov_{q,y}] = \EXP{(\theta-0)(\theta^2 - \vsigma^2)} = 
      \EXP{\theta^3 - \theta \sigma^2 } = 0,\]
and hence $K = [\cov_{q,z}] [\cov_{z,z}]^{-1} = 0,$ and the LBU \refeq{eq:iVII}
will \emph{not} change anything; $q_{a,1}=q_f$.   Looking for the reason for this,
we observe that in the system \refeq{eq:cond-H} in \refT{n-vers}---or in
\refeq{eq:symbolic}---the right-hand-side (rhs) is $\ipj{q \otimes z^{\vee k}}$ in the
$k$-th equation.  In our case this evaluates to
\begin{multline*}
\ipj{q \otimes z^{\vee k}}=\EXP{q \otimes z^{\vee k}} = \EXP{q (q^2+\vepsilon)^k} =\\
  \EXP{q\left(\sum_{i=0}^k c_i q^{2i}\vepsilon^{k-i} \right)} 
 = \sum_{i=0}^k c_i \EXP{\theta^{2i+1}}\EXP{\vepsilon^{k-i}} =0,
\end{multline*}
as $\EXP{\theta^{2i+1}}=0$ for any $i\in\D{N}_0$.  Obviously the $c_i$ are the
binomial coefficients.  This means that in \refeq{eq:cond-H} or \refeq{eq:symbolic}
the rhs vanishes \emph{identically} for any $n\in \D{N}$, and hence
all $\Hf{k}{H}$ will vanish too; i.e.\  no matter what polynomial update we
take, always $\psi_n \equiv 0$, and hence $q_{a,n}=q_f$ for all $n\in \D{N}$.
The loss of information about the sign is so intertwined with the measurement
that no update of the form \refeq{eq:n-deg-pol} can undo it!

If we now come
back to the first idea of choosing a map $\vTheta$ as sketched above, we might chose
$p = |q|$, then $\tilde{Y}(p) = p^2 = |q|^2 = Y(q)$; 
which means we do not care about the sign,
as information about the sign is lost anyway.  The rhs now is---again
neglecting measurement error for the sake of simplicity
\begin{equation*}
\EXP{p \otimes z^{\vee k}} = \EXP{p^{2k+1}} = \EXP{|q|^{2k+1}} =
\EXP{|\theta|^{2k+1}} = \vsigma^{2k+1} 2^k k! \, \sqrt{\frac{2}{\pi}},
\end{equation*}
as these are simply the moments of the \emph{half-normal} or $\chi$-distribution,
and hence one could now compute a polynomial update map $\psi_n$ for any $n$.

One might think that in the formula for the Bayesian update of densities 
\refeq{eq:iIIa} this kind of problem does not appear, but the difficulty comes
when one has
to compute the \emph{likelihood} $p(z|q)$ in \refeq{eq:iIIa}.  Given a measurement
$z = y + \vepsilon$ we have to find \emph{all} $q$ which might have produced it,
and this means that one has to compute the set $Y^{-1}(y)$; so this is where the
difficulty appears then!

We now turn to some examples where we identify parameters in models of varying
complexity.
In \refS{bayes-lin} we will show several examples for the case of $n=1$ for
the update map $\psi_n$, and in \refS{bayes-non-lin} an example for the case
$n=2$.

%
%
%
%
%
%
%
%
%


%

\section{Numerical realisation} \label{S:num-real}
In the instances where we want to employ the theory detailed in the
previous \refS{bayes}, the spaces $\C{U}$ and $\C{Q}$ are usually infinite
dimensional, as is the space $\C{S} = L_2(\vOmega)$.
For an actual computation they have to be discretised or approximated
by finite dimensional spaces.  In our examples we will chose finite
element discretisations and corresponding subspaces.  Hence let
$\C{Q}_M := \text{span }\{\vrho_m\ : m=1,\dots,M\} \subset \C{Q}$ be an
$M$-dimensional subspace with basis $\{\vrho_m\}_{m=1}^M$.  An element
of $\C{Q}_M$ will be represented by the vector $\vek{q}=[q^1, 
\dots, q^M]^T \in \D{R}^M$ such that $\sum^M_{m=1} 
q^m \vrho_m \in \C{Q}_M$.  To avoid a profusion of notations, the corresponding
random vector in $\D{R}^M \otimes \C{S}$ will also be denoted by $\vek{q}$.
The norm $\nd{\vek{q}}_{M}$ to take on $\D{R}^M$
results from the inner product $\bkt{\vek{q_1}}{\vek{q_2}}_{M} := 
\vek{q_1}^T \vek{Q}\vek{q_2}$ with $\vek{Q} = \left(\bkt{\vrho_m}{\vrho_n}_{\C{Q}}\right)$,
the Gram matrix of the basis.  We  will later choose an orthonormal basis, so that
$\vek{Q} = \vek{I}$ is the identity matrix.  Similarly, on 
$\E{Q}_M = \D{R}^M \otimes \C{S}$ the inner product is 
$\bkt{\vek{q}_1}{\vek{q}_2}_{\E{Q}_M}:= \EXP{\bkt{\vek{q_1}}{\vek{q_2}}_{M}}$.
The space of possible measurements
can usually be taken to be finite dimensional (here $=R$), whose elements are
similarly represented by a vector of coefficients $\vek{z} \in \D{R}^R$.

On $\D{R}^M$, representing $\C{Q}_M$, the Kalman gain operator in \refT{n-1-vers}
in \refeq{eq:iVII} becomes a matrix $\vek{K}\in \D{R}^{M \times R}$.
Then the update corresponding to \refeq{eq:iVII} is
\begin{equation}  \label{eq:iIX}
  \vek{q}_a = \vek{q}_f + \vek{K}(\vhat{z} - \vek{z}), \text{   with   } 
  \vek{K} = \vek{C}_{q,z}\, \vek{C}_{z,z}^{-1}.
\end{equation}
Here the covariances are $\vek{C}_{q,z} := \EXP{\tilde{\vek{q}}\;
 \tilde{\vek{z}}^T} = \EXP{\tilde{\vek{q}} \otimes \tilde{\vek{z}}}$,
and similarly for $\vek{C}_{z,z}$.  Often the measurement error $\vek{\vepsilon}$ is
independent of $\vek{q}$ --- actually uncorrelated would be \emph{sufficient}---hence
$\vek{C}_{z,z} = \vek{C}_{y,y} + \vek{C}_{\vepsilon,\vepsilon}$ and $\vek{C}_{q,z} =
\vek{C}_{q,y}$.  We once more
recall our comments in \refSS{approx-cond-expect} following \refeq{eq:LBU-n1}
regarding the inverse which also appears in \refeq{eq:iIX}.  Recall that usually
the error model involves a regular covariance $\vek{C}_{\vepsilon,\vepsilon}$,
so that $\vek{C}_{z,z} = \vek{C}_{y,y} + \vek{C}_{\vepsilon,\vepsilon}$ is at
least theoretically regular.

It is important to emphasise that the theory presented in the forgoing \refS{bayes}
is independent of any discretisation.  But one usually can still 
not numerically compute with objects like $\vek{q}\in
\E{Q}_M = \D{R}^M \otimes \C{S}$, as $\C{Q} = L_2(\vOmega)$ is normally an
infinite dimensional space and has to be discretised.  One well-known possibility
are samples, i.e.\ the RV $\vek{q}(\omega)$ is represented by its value at certain points
$\omega_z$, and the points usually come from some quadrature rule.  The well-known
Monte Carlo (MC) method uses random samples, the quasi-Monte Carlo (QMC) method uses
low discrepancy samples, and other rules like sparse grids (Smolyak rule) are possible.
Using MC samples in the context of the linear update \refeq{eq:iVII}
is known as the \emph{Ensemble Kalman Filter} (EnKF), see \citep{BvrAkJsOpHgm11}
for a general overview in this context, and \citep{Evensen2009, Evensen2009a} for
a thorough description and analysis.  This method is concepyually fairly simple
and is currently 
a favourite for problems where the computation of the predicted measurement
$\vek{y}_f(\omega_z)$ is difficult or expensive.  It needs far fewer samples
for meaningful results than MCMC, but on the other hand it uses the linear
approximation inherent in \refeq{eq:iIX}.

Here we want to use so-called \emph{functional} or \emph{spectral} approximations, so
similarly as for $\C{Q}_M$, we pick a finite set of linearly independent
vectors in $\C{S}$.  As $\C{S} = L_2(\vOmega)$, these abstract vectors
are in fact RVs with finite variance.  Here we will use the best known
example, namely \emph{Wiener}'s \emph{polynomial chaos} expansion (PCE) as basis 
\citep{Wiener1938, ghanemSpanos91, holdenEtAl96, Janson1997, malliavin97, matthies6},
this allows us to use \refeq{eq:iIX} without sampling, see 
\citep{BvrAkJsOpHgm11, opBvrAlHgm12, bvrAlOpHgm12-a, boulder:2011, OpBrHgm12}, and also
\citep{saadGhn:2009, Blanchard2010a}.

The PCE is an expansion in multivariate \emph{Hermite polynomials} 
\citep{ghanemSpanos91, holdenEtAl96, Janson1997, malliavin97, matthies6}; we denote
by $H_{\vek{\alpha}}(\vek{\theta}) = \prod_{k \in \D{N}} h_{\alpha_k}(\theta_k)
\in \C{S}$ the multivariate polynomial in standard independent Gaussian RVs 
$\vek{\theta}(\omega) = (\theta_1(\omega),\dots, \theta_k(\omega), 
\dots)_{k\in \D{N}}$,
where $h_j$ is the usual univariate Hermite polynomial, and $\vek{\alpha} = 
(\alpha_1,\dots,\alpha_k,\dots)_{k\in \D{N}}\in \C{N}:=\D{N}_0^{(\D{N})}$
is a multi-index of generally infinite lenght but with only finitely many
entries non-zero.  As $h_0 \equiv 1$,
the infinite product is effectively finite and always well-defined.

The \emph{Cameron-Martin} theorem assures us \citep{holdenEtAl96, malliavin97, Janson1997}
that the set of these polynomials is dense in $\C{S} = L_2(\vOmega)$,
and in fact $\{H_{\vek{\alpha}}/\sqrt{(\vek{\alpha} !)} \}_{\vek{\alpha} 
\in \C{N}}$ is a complete orthonormal system (CONS), 
where $\vek{\alpha} ! := \prod_{k \in \D{N}} (\alpha_k !)$ is the product
of the individual factorials, also well-defined as except for finitely many 
$k$ one has $\alpha_k ! = 0! = 1$.  So we may write $\vek{q}(\omega) =
\sum_{\vek{\alpha}\in \C{N}} \vek{q}^{\vek{\alpha}} H_{\vek{\alpha}}
(\vek{\theta}(\omega))$ with $\vek{q}^{\vek{\alpha}} \in \D{R}^M$,
and similarly for $\vek{z}$ and all other RVs.
In this way the RVs are expressed as functions of other,
known RVs $\vek{\theta}$---hence the name \emph{functional} 
approximation---and not through samples.

The space $\C{S}$ may now be discretised by taking a finite subset $\C{J}
\subset \C{N}$ of size $J = \ns{\C{J}}$, and setting $\C{S}_J = \text{span }
\{H_{\vek{\alpha}}\,:\, \vek{\alpha} \in \C{J} \} \subset \C{S}$.  The
orthogonal projection $P_J$ onto $\C{S}_J$ is then simply
\begin{equation}  \label{eq:proj-J}
P_J: \C{Q}_M \otimes \C{S} \ni
\sum_{\vek{\alpha}\in \C{N}} \vek{q}^{\vek{\alpha}} H_{\vek{\alpha}} \mapsto
\sum_{\vek{\alpha}\in \C{J}} \vek{q}^{\vek{\alpha}} H_{\vek{\alpha}} 
\in \C{Q}_M \otimes \C{S}_J.
\end{equation}
We then take \refeq{eq:iIX} and rewrite it as
\begin{eqnarray}  \label{eq:proj-lin-f1}
  \vek{q}_a &=& \vek{q}_f + \vek{K}(\vek{z} - \vek{y}_f) =\\
  \sum_{\vek{\alpha}\in \C{N}} \vek{q}_a^{\vek{\alpha}} H_{\vek{\alpha}}(\vek{\theta}) &=& 
  \sum_{\vek{\alpha}\in \C{N}} \left(\vek{q}_f^{\vek{\alpha}} + \vek{K}\left( 
  \vek{z}^{\vek{\alpha}}-\vek{y}^{\vek{\alpha}}_f\right)\right)
  H_{\vek{\alpha}}(\vek{\theta}). \label{eq:proj-lin-f2}
\end{eqnarray}
Projecting both sides of \refeq{eq:proj-lin-f2} is very
simple and results in
\begin{equation} \label{eq:proj-lin-J}
  \sum_{\vek{\alpha}\in \C{J}} \vek{q}_a^{\vek{\alpha}} H_{\vek{\alpha}} = 
  \sum_{\vek{\alpha}\in \C{J}} \left(\vek{q}_f^{\vek{\alpha}} + \vek{K}\left( 
  \vek{z}^{\vek{\alpha}}-\vek{y}^{\vek{\alpha}}_f\right)\right)H_{\vek{\alpha}}.
\end{equation}
Obviously the projection $P_J$ commutes with the Kalman operator $K$ and
hence with its finite dimensional analogue $\vek{K}$.  One may actually
concisely write \refeq{eq:proj-lin-J} as
\begin{equation} \label{eq:proj-comm-K}
  P_J \vek{q}_a = P_J \vek{q}_f + P_J \vek{K}(\vek{z} - \vek{y}_f) =
  P_J\vek{q}_f + \vek{K}(P_J\vek{z} - P_J\vek{y}_f).
\end{equation}

Elements of the discretised space $\E{Q}_{M,J} = \C{Q}_M \otimes
\C{S}_J \subset \E{Q}$ thus may be written as 
$\sum_{m=1}^M \sum_{\vek{\alpha}\in \C{J}} q^{\vek{\alpha},m} \vrho_m 
H_{\vek{\alpha}}$.  The tensor representation is
$\mat{q} := (q^{\vek{\alpha},m}) = \sum_{\vek{\alpha}\in \C{J}} 
\vek{q}^{\vek{\alpha}} \otimes \vek{e}^{\vek{\alpha}}$, where the 
$\vek{e}^{\vek{\alpha}}$ are the unit vectors in $\D{R}^J$, may be used
to express \refeq{eq:proj-lin-J} or \refeq{eq:proj-comm-K} succinctly as
\begin{equation}  \label{eq:proj-t}
 \mat{q}_a = \mat{q}_f+ \mat{K}(\mat{z}-\mat{y}_f),
\end{equation}
again an equation between the tensor representations of some RVs,
where $\mat{K} = \vek{K} \otimes \vek{I}$ with $\vek{K}$ from \refeq{eq:iIX}.
Hence the update equation is naturally in a tensorised form.  This is
how the update can finally be computed in the PCE representation without
any sampling \citep{BvrAkJsOpHgm11, opBvrAlHgm12, bvrAlOpHgm11, boulder:2011}.
Analogous statements hold for the forms of the update \refeq{eq:n-deg-pol}
with higher order terms $n>1$, and do not have to be repeated here.  Let us
remark that these updates go very seamlessly with very efficient methods for
sparse or low-rank approximation of tensors, c.f.\ the monograph
\citep{Hackbusch_tensor} and the literature therein.  These methods are PCE-forms
of the Bayesian update, and in particular the \refeq{eq:proj-t}, because of
its formal affinity to the Kalman filter (KF), may be called the polynomial
chaos expansion based Kalman filter (PCEKF).

It remains to say how to compute the terms $\Hf{k}{H}$ in the update equation 
\refeq{eq:n-deg-pol}---or rather the terms in the defining \refeq{eq:cond-H}
in \refT{n-vers}---in this approach.  Given the PCEs of the RVs, this is actually
quite simple as any moment can be computed directly from the PCE
\citep{matthies6, opBvrAlHgm12, bvrAlOpHgm12-a}. 
A typical term $\ipj{z^{\vee k}} = \ipj{\text{Sym}(z^{\otimes k})}=
\text{Sym}(\ipj{z^{\otimes k}})$
in the operator matrix \refeq{eq:cond-H}, where $\vek{z}=\sum_\alpha\vek{z}^{\vek{\alpha}}
H_{\vek{\alpha}}(\vek{\theta})$, may be computed through
\begin{multline}  \label{eq:typ-zk}
  \ipj{\vek{z}^{\otimes k}} = \EXP{\bigotimes_{i=1}^k \sum_{\vek{\alpha}_i} \left( 
  \vek{z}^{\vek{\alpha}_i} H_{\vek{\alpha}_i}\right)} = \\
   \EXP{\sum_{\vek{\alpha}_1, \dots, \vek{\alpha}_k} \bigotimes_{i=1}^k 
   \vek{z}^{\vek{\alpha}_i} \prod_{i=1}^k H_{\vek{\alpha}_i}}  = 
   \sum_{\vek{\alpha}_1, \dots, \vek{\alpha}_k} \bigotimes_{i=1}^k \vek{z}^{\vek{\alpha}_i}
   \;\EXP{\prod_{i=1}^k H_{\vek{\alpha}_i}}
\end{multline}
As here the $H_{\vek{\alpha}}$ are \emph{polynomials}, 
the last expectation in \refeq{eq:typ-zk}
is finally over products of powers of pairwise independent normalised Gaussian variables,
which actually may be done analytically \citep{holdenEtAl96, malliavin97, Janson1997}.
But some simplifications come from
remembering that $\vek{z}^0=\EXP{\vek{z}} = \bar{\vek{z}}$, $H_{\vek{0}}\equiv 1$, 
the orthogonality relation $\bkt{H_{\vek{\alpha}}}{H_{\vek{\beta}}} = 
\delta_{\vek{\alpha},\vek{\beta}}\, \vek{\alpha}!$, and that the Hermite polynomials
are an \emph{algebra}.  Hence $H_{\vek{\alpha}}H_{\vek{\beta}} = \sum_{\vek{\gamma}}
c^{\vek{\gamma}}_{\vek{\alpha},\vek{\beta}}H_{\vek{\gamma}}$, where the \emph{structure}
coefficients $c^{\vek{\gamma}}_{\vek{\alpha},\vek{\beta}}$ are known analytically
\citep{malliavin97, matthies6, opBvrAlHgm12, bvrAlOpHgm12-a}.

Similarly, for a typical right-hand-side term $\ipj{q\otimes z^{\vee k}} = \ipj{q\otimes
\text{Sym}(z^{\otimes k})}$ in \refeq{eq:cond-H} with $\vek{q}=\sum_\beta 
\vek{q}^{\vek{\beta}} H_{\vek{\beta}}(\vek{\theta})$ one has
\begin{equation}  \label{eq:typ-qzk}
   \ipj{q\otimes \text{Sym}(z^{\otimes k})} =
     \sum_{\vek{\beta}, \vek{\alpha}_1, \dots, \vek{\alpha}_k} \vek{q} \otimes
     \text{Sym}\left(\bigotimes_{i=1}^k \vek{z}^{\vek{\alpha}_i}\right)
   \;\EXP{H_{\vek{\beta}} \, \prod_{i=1}^k H_{\vek{\alpha}_i}}.
\end{equation}
As these relations may seem a bit involved---they are actually just a bit intricate
combination of \emph{known} terms---we show here how simple they become for
the case of the covariance needed in the linear update formula \refeq{eq:iVII}
or rather \refeq{eq:iIX}:
\begin{eqnarray}  \label{eq:cov-PCE-1}
   \vek{C}_{z,z} &=
   \sum_{\vek{\alpha}\in \C{N}, \vek{\alpha} \ne 0} (\vek{\alpha} !)\;
    \vek{q}^{\vek{\alpha}}\otimes \vek{z}^{\vek{\alpha}} 
   &\approx 
   \sum_{\vek{\alpha}\in \C{J}, \vek{\alpha} \ne 0} (\vek{\alpha} !)\;
   \vek{z}^{\vek{\alpha}}\otimes \vek{z}^{\vek{\alpha}},\\
  \vek{C}_{q,z} &= 
   \sum_{\vek{\alpha}\in \C{N}, \vek{\alpha} \ne 0} (\vek{\alpha} !)\;
    \vek{q}^{\vek{\alpha}}\otimes \vek{z}^{\vek{\alpha}} 
  &\approx \sum_{\vek{\alpha}\in \C{J}, \vek{\alpha} \ne 0} (\vek{\alpha} !)\;
   \vek{q}^{\vek{\alpha}}\otimes \vek{z}^{\vek{\alpha}}. 
    \label{eq:cov-PCE-2}
\end{eqnarray}

Looking for example at \refeq{eq:iIX} and our setup as explained in \refS{intro},
we see that the coefficients of $\vek{z}=\sum_\alpha\vek{z}^{\vek{\alpha}}
H_{\vek{\alpha}}$ or rather those of $\vek{y}=\sum_\alpha\vek{y}^{\vek{\alpha}}
H_{\vek{\alpha}} = \vek{Y}(\vek{q})$ have to be computed from those of
$\vek{q}=\sum_\beta\vek{q}^{\vek{\beta}} H_{\vek{\beta}}$.  This propagation
of uncertainty through the system is known as \emph{uncertainty quantification} (UQ),
e.g.\ \citep{matthies6} and the references therein.  For the sake of brevity,
we will not touch further on this subject, which nevertheless is the bedrock
on which we built the whole computational procedure.

We next concentrate in \refS{bayes-lin} on examples of updating with $\psi_n$ 
for the case $n=1$ in \refeq{eq:n-deg-pol}, whereas in \refS{bayes-non-lin}
an example for the case $n=2$ in \refeq{eq:n-deg-pol} will be shown.




%
%
%
%
%


%

\section{The linear Bayesian update} \label{S:bayes-lin}
All the examples in this section have been computed with the case $n=1$
of up to linear terms in \refeq{eq:n-deg-pol}, i.e.\ this is the LBU with
PCEKF.
\begin{figure}[!ht]
\centering
 \includegraphics[width=0.6\textwidth,height=0.35\textheight]{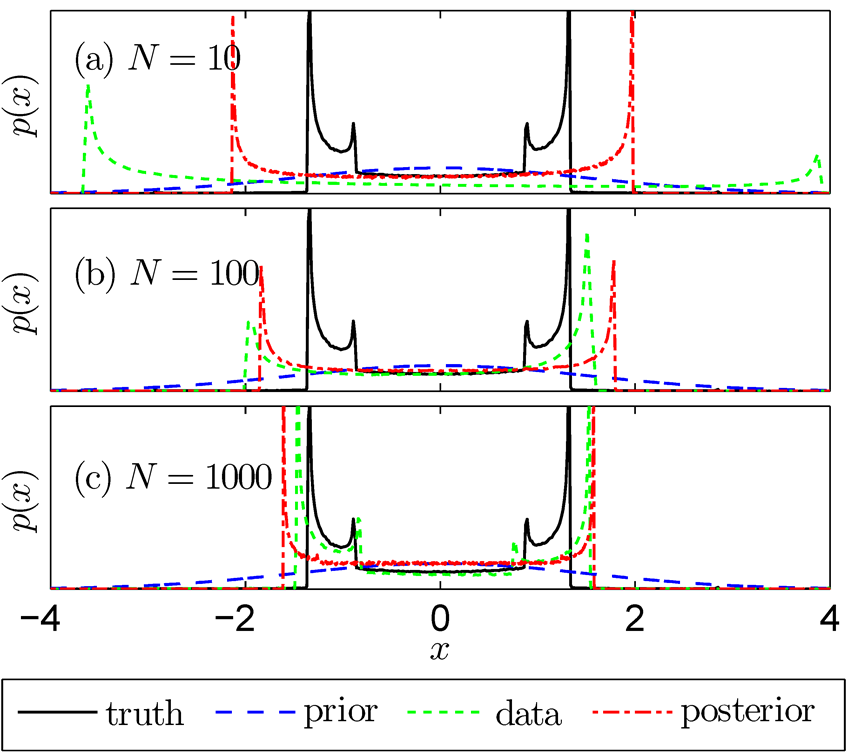}
 \caption{pdfs for linear Bayesian update \cite{opBvrAlHgm12}}
\label{F:exp-B-1}
\end{figure}
As the traditional Kalman filter is highly geared towards Gaussian
distributions \citep{Kalman}, and also its Monte Carlo variant EnKF
which was mentioned in \refS{num-real} tilts towards Gaussianity,
we start with a case---already described in \citep{opBvrAlHgm12}---where
the the quantity to be identified has a strongly
non-Gaussian distribution, shown in black---the `truth'---in \refig{exp-B-1}.
The operator describing the system is the identity---we compute the quantity
directly, but there is a Gaussian measurement error.  The `truth' was
represented as a $12^{\text{th}}$ degree PCE.
We use the methods as described in \refS{num-real}, and here in particular
the \refeq{eq:iIX} and \refeq{eq:proj-t}, the PCEKF.

The update is repeated several times (here ten times)
with new measurements---see \refig{exp-B-1}.  The task is here
to identify the distribution labelled as `truth' with ten updates of $N$
samples (where $N=10, 100, 1000$ was used), and we start with
a very broad Gaussian prior (in blue).  Here we see the ability of
the polynomial based LBU, the PCEKF, to identify highly non-Gaussian distributions,
the posterior is shown in red and the pdf estimated from the samples in green;
for further details see \citep{opBvrAlHgm12}.

\begin{figure}[!ht]
\centering
 \includegraphics[width=0.9\textwidth,height=0.35\textheight]{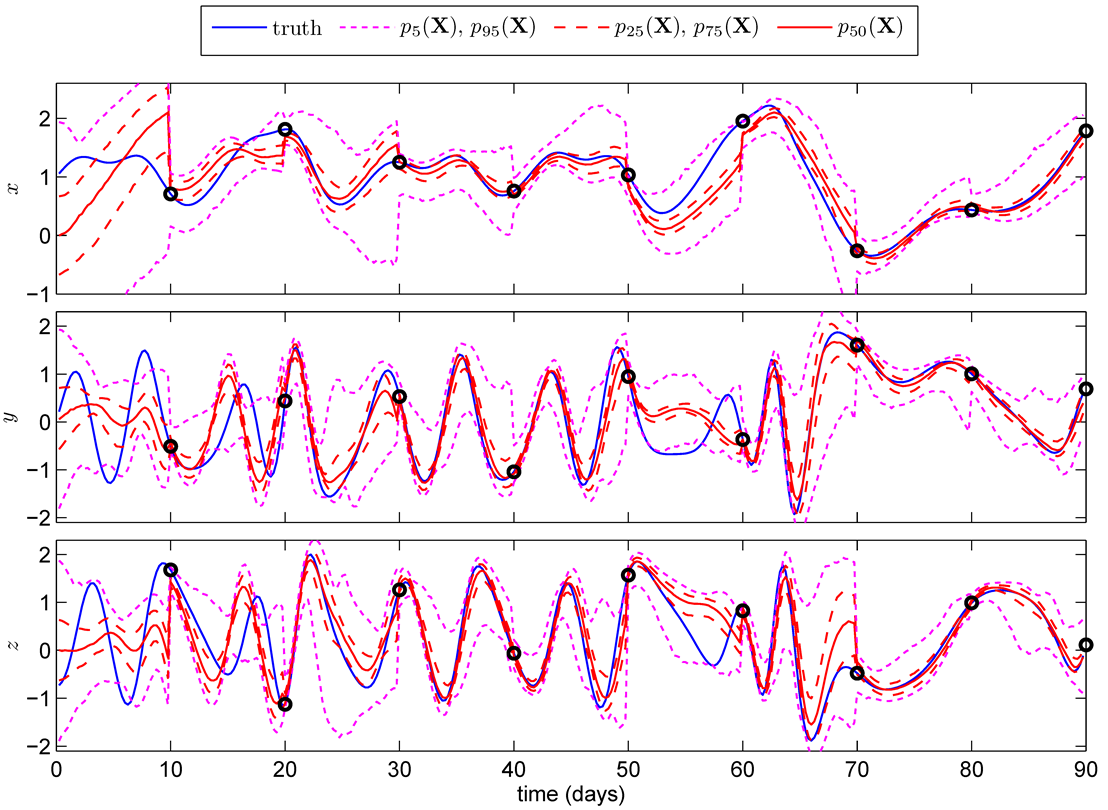}
 \caption{Time evolution of Lorenz-84 state and uncertainty with LBU \cite{opBvrAlHgm12}}
\label{F:exp-B-2}
\end{figure}
The next example is also from \citep{opBvrAlHgm12}, where the system is the well-known
Lorenz-84 chaotic model, a system of three nonlinear ordinary differential equations
operating in the chaotic regime.  Remember that this was originally a model to
describe the evolution of some amplitudes of a spherical harmonic expansion of
variables describing world climate.  As the original scaling of the variables has
been kept, the time axis in \refig{exp-B-2} is in \emph{days}.  Every ten days a
noisy measurement is performed and the state description is updated.  In between
the state description evolves according to the chaotic dynamic of the system.
One may observe from \refig{exp-B-2} how the uncertainty---the width of the distribution
as given by the quantile lines---shrinks every time a measurement is performed, and
then increases again due to the chaotic and hence noisy dynamics.
Of course, we did not really measure world climate, but rather simulated the `truth'
as well, i.e.\ a \emph{virtual} experiment, like the others to follow.
More details may be found in \citep{opBvrAlHgm12} and the references therein.

\begin{figure}[!ht]
\centering
\begin{minipage}{.4\textwidth}
  \centering
  \includegraphics[width=.99\linewidth]{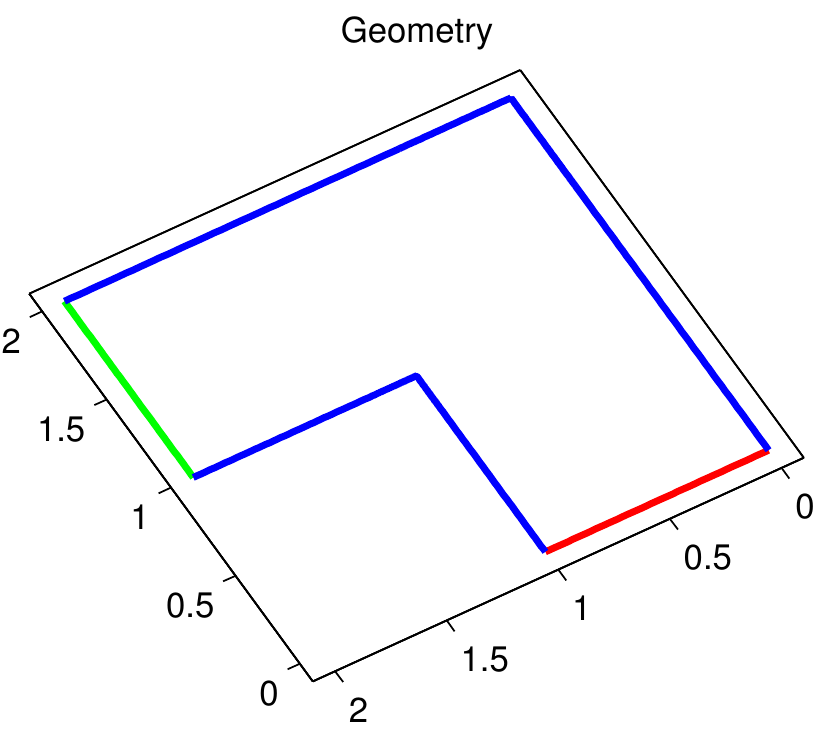}
  \captionof{figure}{Diffusion domain}
  \label{F:exp-B-3}
\end{minipage}%
\begin{minipage}{.6\textwidth}
  \centering
  \includegraphics[width=.99\linewidth]{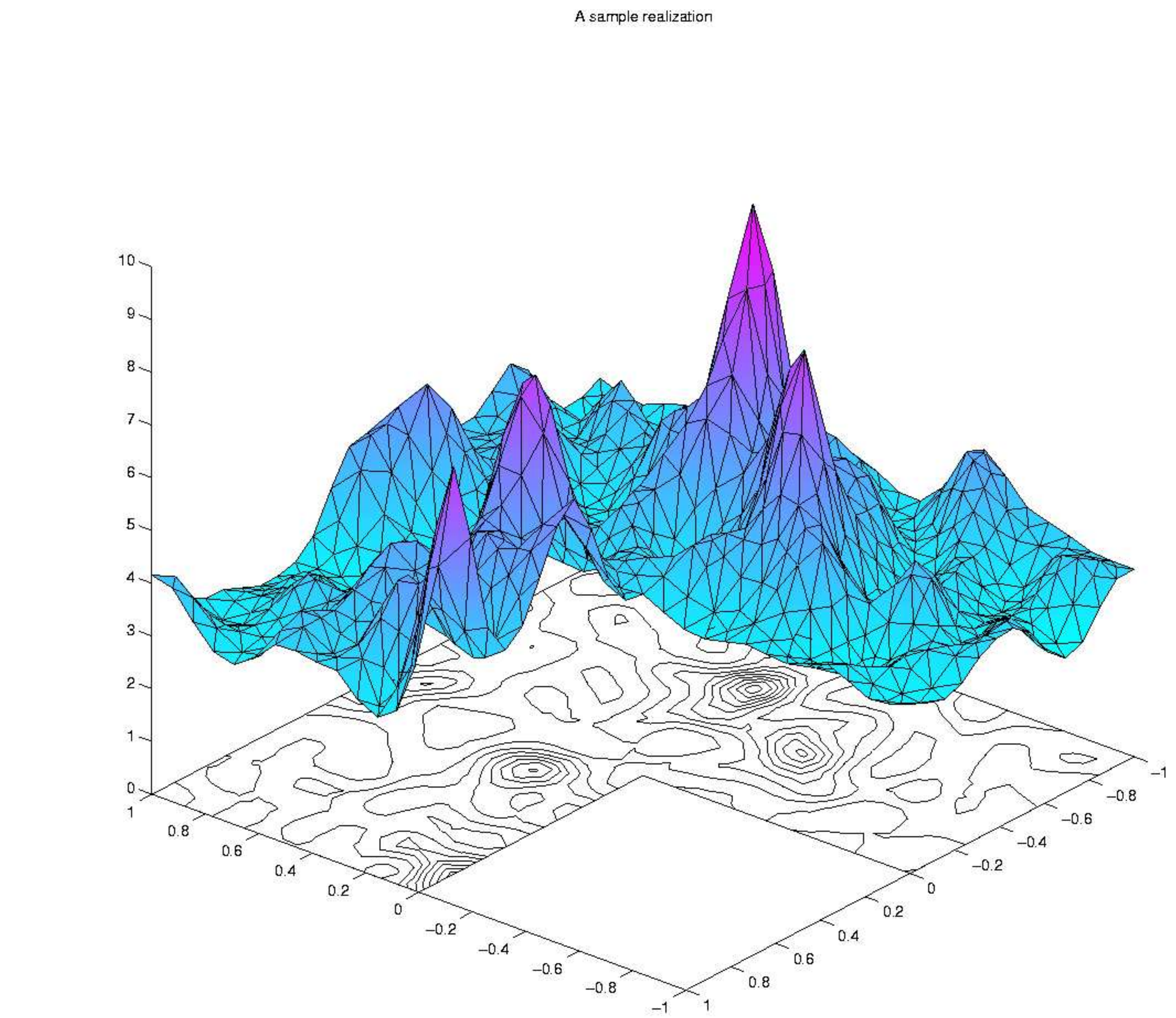}
  \captionof{figure}{Conductivity field}  
  \label{F:exp-B-4}
\end{minipage}
\end{figure}
From \citep{bvrAlOpHgm11, bvrAlOpHgm12-a} we take the example shown in
\refig{exp-B-3}, a linear stationary diffusion equation on an L-shaped
plane domain as alluded to in \refS{intro}.  The diffusion coefficient
$\kappa$ in \refeq{eq:I-c} is to be identified.  As argued in \citep{BvrAkJsOpHgm11},
it is better to work with $q = \log \kappa$ as the diffusion coefficient has
to be positive, but the results are shown
in terms of $\kappa$.

\begin{figure}[!ht]
\centering
\begin{minipage}{0.5\textwidth}
  \centering
  \includegraphics[width=0.7\linewidth]{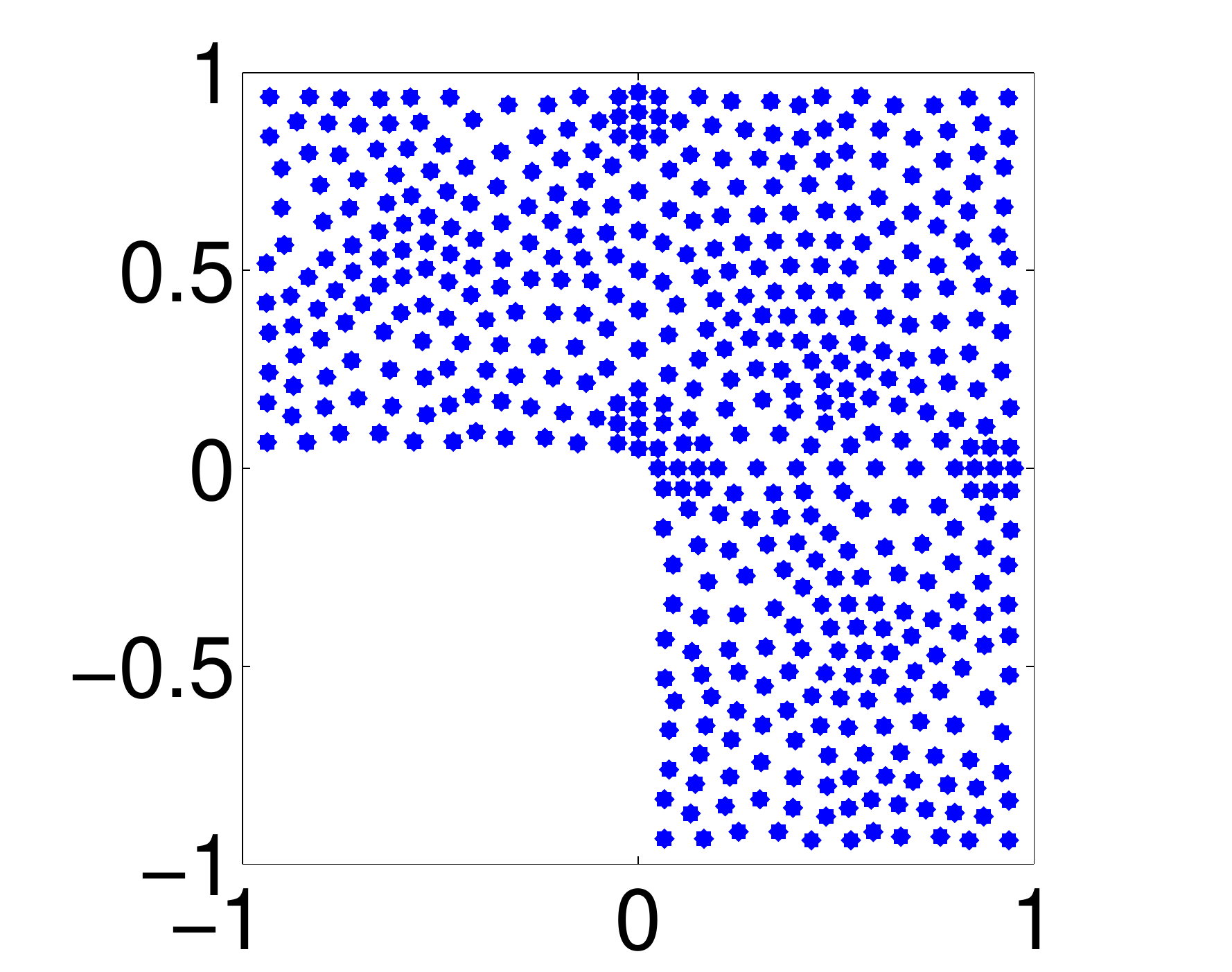}
  \captionof{figure}{447 measurement patches}
  \label{F:exp-B-5}
\end{minipage}%
\begin{minipage}{0.5\textwidth}
  \centering
  \includegraphics[width=0.7\linewidth]{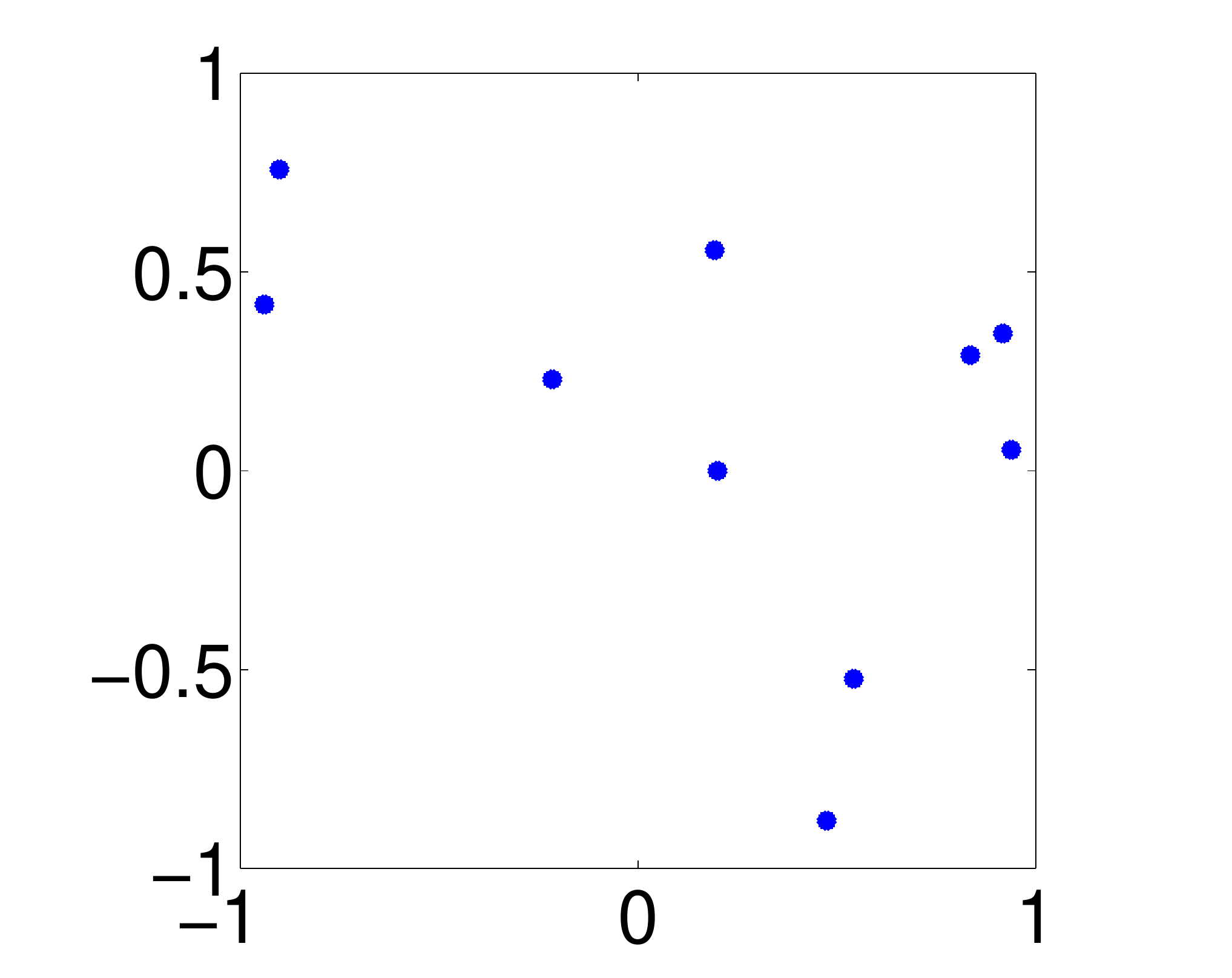}
  \captionof{figure}{10 measurement patches}  
  \label{F:exp-B-6}
\end{minipage}
\end{figure}
One possible realisation of the diffusion coefficient
is shown in \refig{exp-B-4}.  More realistically, one should assume that
$\kappa$ is a symmetric positive definite tensor field, unless one knows that
the diffusion is \emph{isotropic}.  Also in this case one should do the updating
on the logarithm.  For the sake of simplicity we stay with the scalar case,
as there is no principal novelty in the non-isotropic case.

The virtual experiments use different right-hand-sides $f$ in \refeq{eq:I-c},
and the measurement is the observation of the solution $u$ averaged over little patches,
two of these arrangements are shown in \refig{exp-B-5} and \refig{exp-B-6}.

\begin{figure}[!ht]
\centering
\begin{minipage}{0.5\textwidth}
  \centering
  \includegraphics[width=0.99\linewidth]{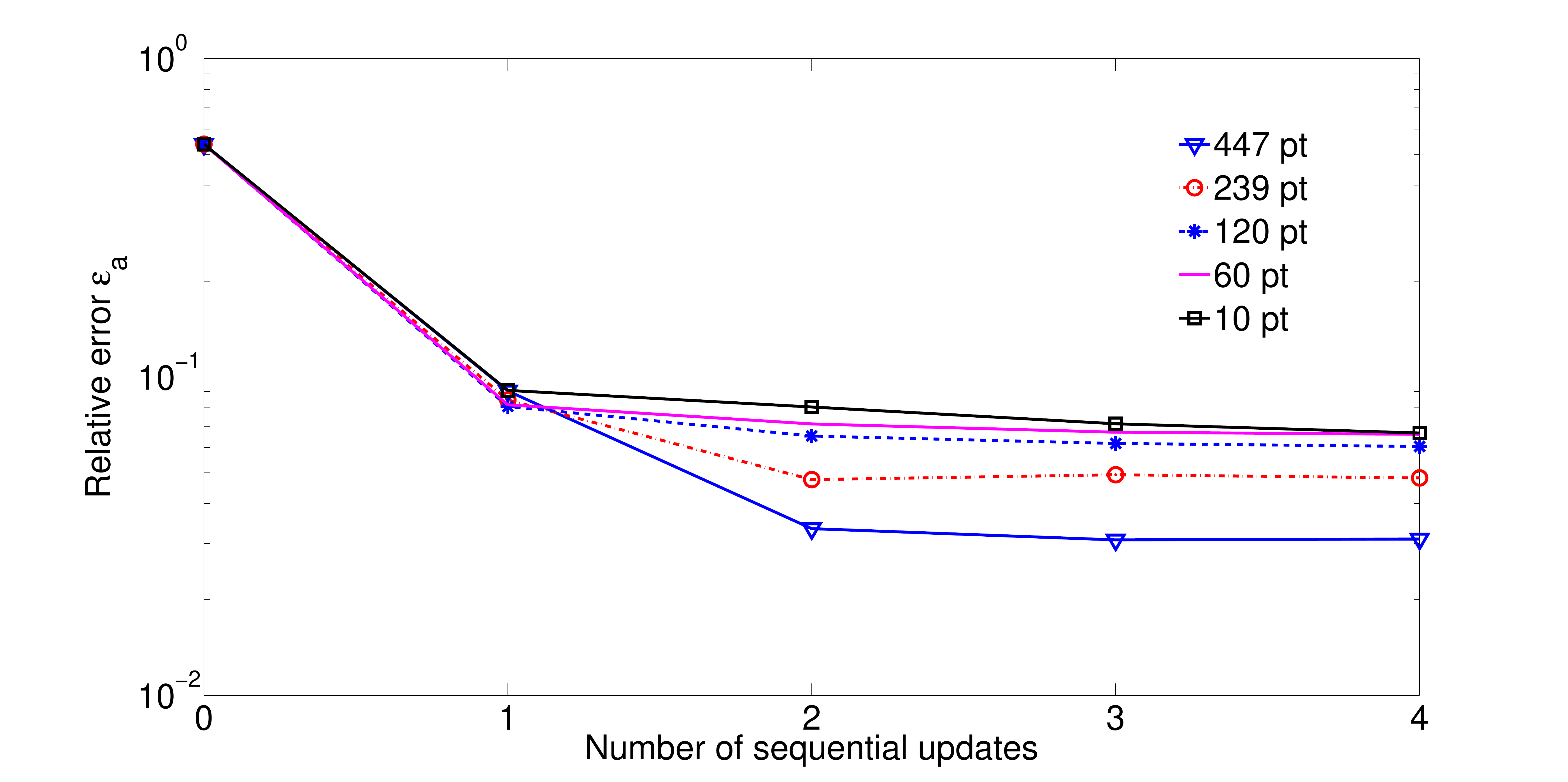}
  \captionof{figure}{Convergence of identification \cite{bvrAlOpHgm12-a}}
  \label{F:exp-B-7}
\end{minipage}%
\begin{minipage}{0.5\textwidth}
  \centering
  \includegraphics[width=0.99\linewidth,height=0.17\textheight]{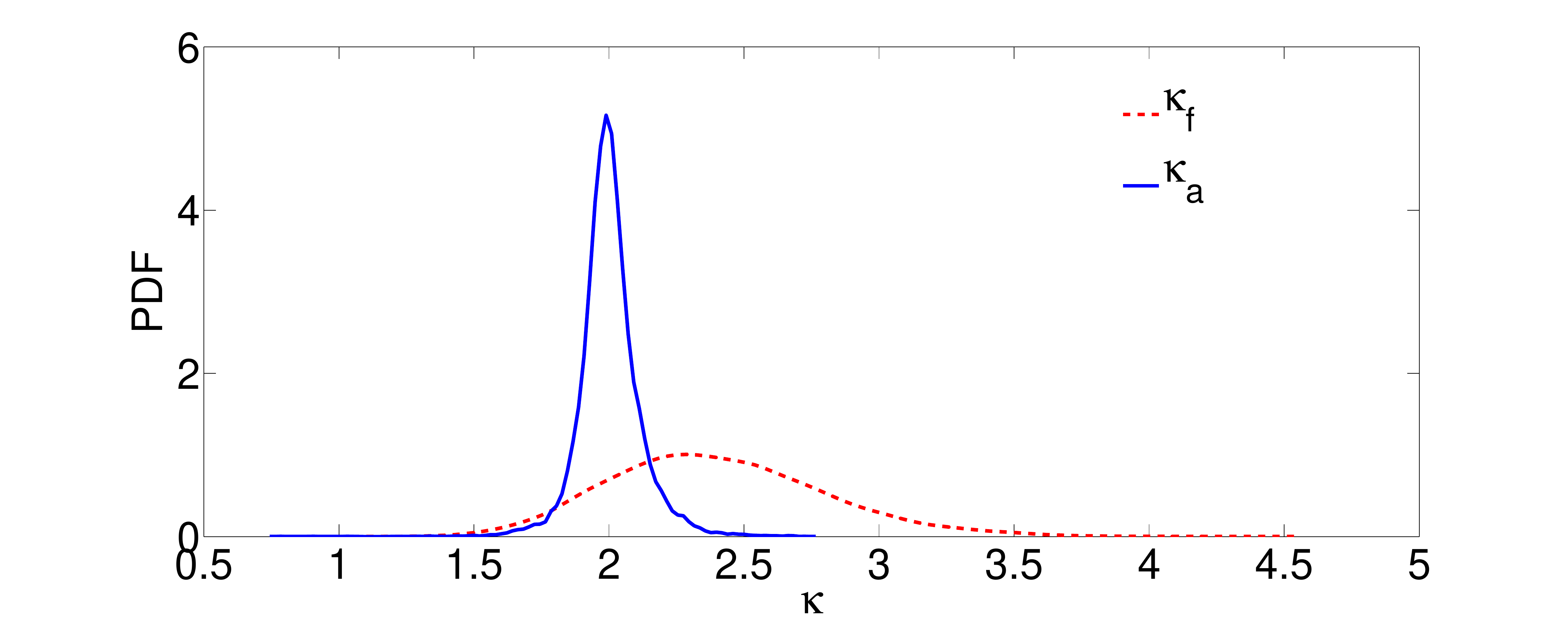}
  \captionof{figure}{Prior and posterior \cite{bvrAlOpHgm12-a}}  
  \label{F:exp-B-8}
\end{minipage}
\end{figure}
In \refig{exp-B-7} one may observe the decrease of the error with successive
updates, but due to measurement error and insufficient information from just
a few patches, the curves level off, leaving some residual uncertainty.
The pdfs of the diffusion coefficient at some point in the domain before
and after the updating is shown in \refig{exp-B-8}, the `true' value at
that point was $\kappa=2$.
Further details can be found in \citep{bvrAlOpHgm11, bvrAlOpHgm12-a}.

\begin{figure}[!ht]
\centering
\begin{minipage}{0.5\textwidth}
  \centering
  \includegraphics[width=0.99\linewidth]{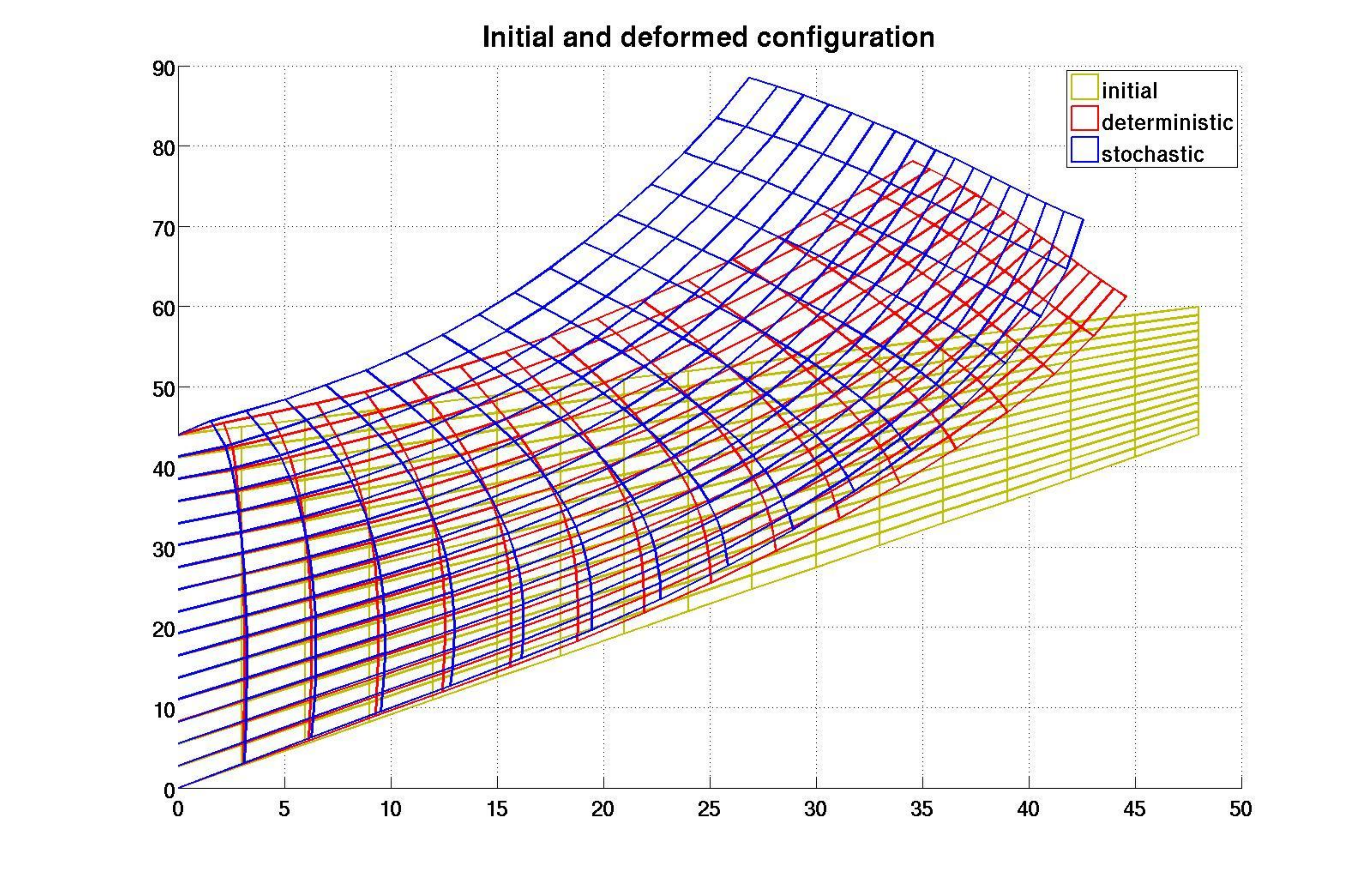}
  \captionof{figure}{Deformations \cite{RosicPhD}}
  \label{F:exp-B-9}
\end{minipage}%
\begin{minipage}{0.5\textwidth}
  \centering
  \includegraphics[width=0.99\linewidth]{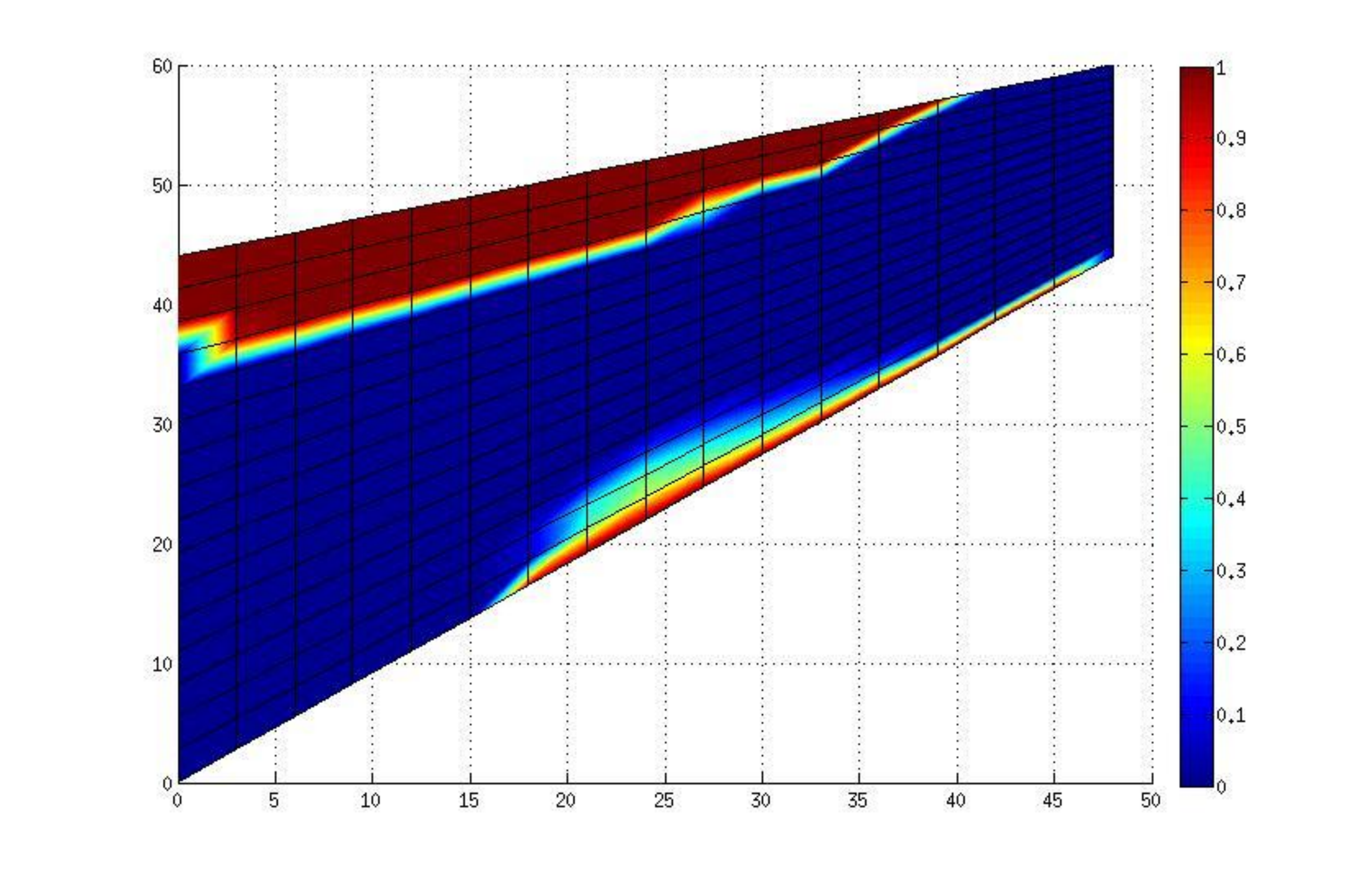}
  \captionof{figure}{Exceedance probability \cite{RosicPhD}}  
  \label{F:exp-B-10}
\end{minipage}
\end{figure}
As a last example with LBU, we take a strongly nonlinear and also non-smooth
situation, namely elasto-plasticity with linear hardening and large deformations
and a \emph{Kirchhoff-St.\ Venant} elastic material law 
\citep{RosicPhD, BvrAkJsOpHgm11, rosic2013hgm}.  This example is known as \emph{Cook's membrane},
and is shown in \refig{exp-B-9} with the undeformed mesh (initial), the deformed one
obtained by computing with average values of the elasticity and plasticity material
constants (deterministic), and finally the average result from a stochastic forward
calculation of the probabilistic model (stochastic), which is described by a
variational inequality \citep{RosicPhD, rosic2013hgm}.  In \refig{exp-B-10} one may get another
impression of results of the forward model, the probability of the \emph{von Mises}
stress being beyond a certain value.

\begin{figure}[!ht]
\centering
 \includegraphics[width=0.9\textwidth]{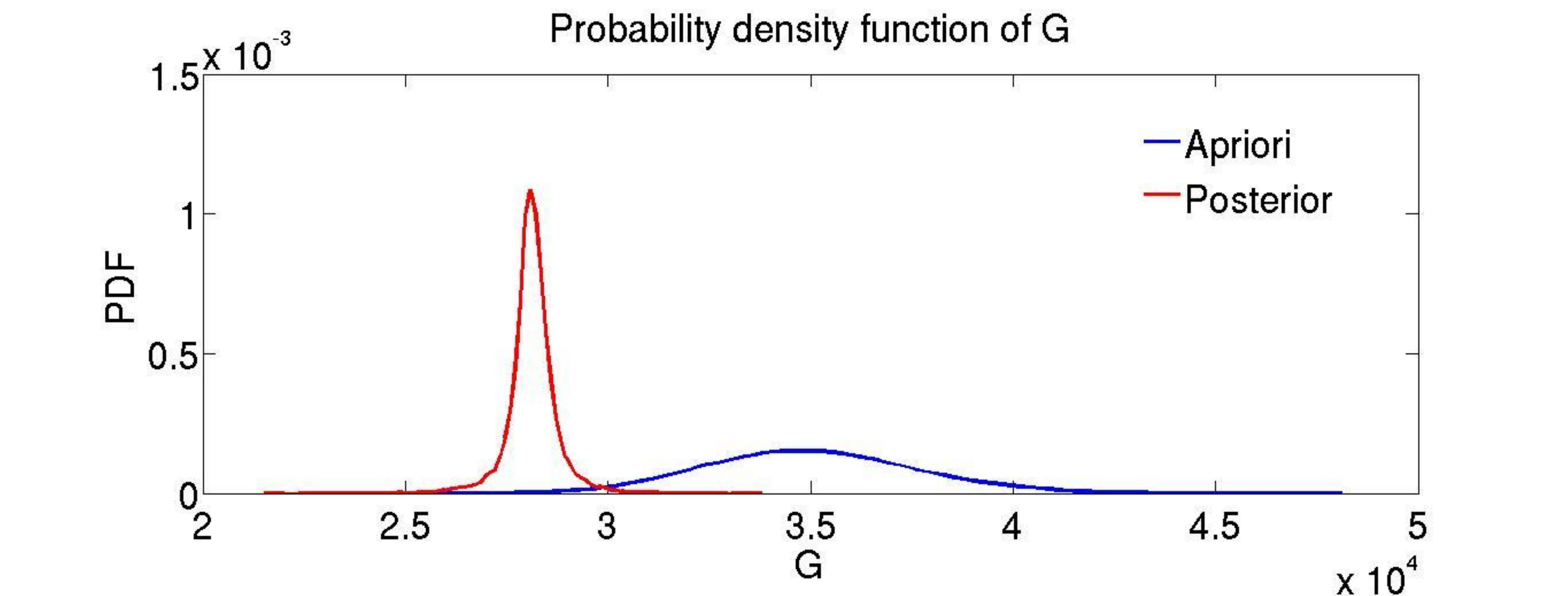}
 \caption{Prior and posterior of shear modulus \cite{rosic2013hgm}}
\label{F:exp-B-11}
\end{figure}
The shear modulus $G$ has to be identified, which is made more difficult by the
non-smooth nonlinearity.  In \refig{exp-B-11} one may see the prior and posterior
distributions of the shear modulus at one point in the domain.  The `truth' is
$G \approx 2.7$, and one may observe that the update is successful although
the prior density almost vanishes at $G=2.7$.

%
%
%
%
%
%
%


%

\section{The nonlinear Bayesian update} \label{S:bayes-non-lin}
\begin{figure}[!ht]
\centering
 \includegraphics[width=0.9\textwidth,height=0.2\textheight]{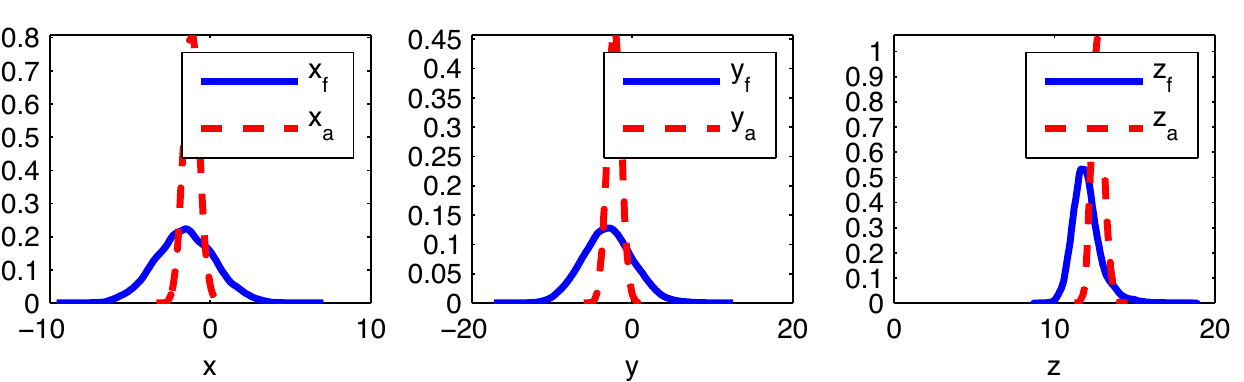}
 \caption{Linear measurement: prior and posterior after one update}
\label{F:exp-NB-5}
\end{figure}
In this Section we want to show a computation with the case $n=2$
of up to quadratic terms in \refeq{eq:n-deg-pol}.  We go back to
the example of the chaotic Lorentz-84 \citep{opBvrAlHgm12} model already shown in
\refS{bayes-lin}.  For this kind of experiment it has several advantages:
it has only a three-dimensional state space, these are the uncertain `parameters',
i.e.\ $(x,y,z)\in\C{Q}=\D{R}^3$, the corresponding operator $A$ in the 
abstract \refeq{eq:I} is sufficiently nonlinear to make the problem difficult,
and adding to this we operate the equation in its chaotic regime, so that
new uncertainty is added between measurements.

\begin{figure}[!ht]
\centering
 \includegraphics[width=0.9\textwidth,height=0.2\textheight]{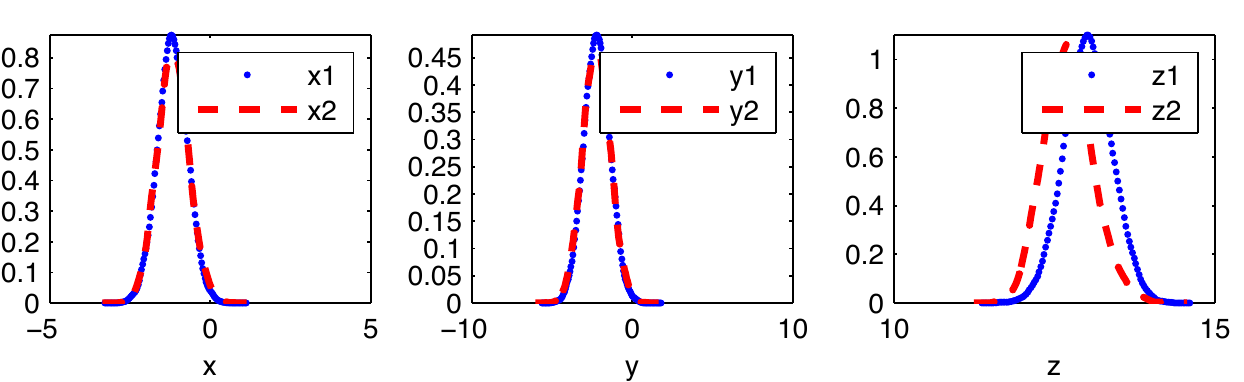}
 \caption{Linear measurement: Comparison posterior for LBU ($n=1$) and 
 NLBU ($n=2$) after one update}
\label{F:exp-NB-5.5}
\end{figure}
\begin{figure}[!ht]
\centering
 \includegraphics[width=0.9\textwidth,height=0.2\textheight]{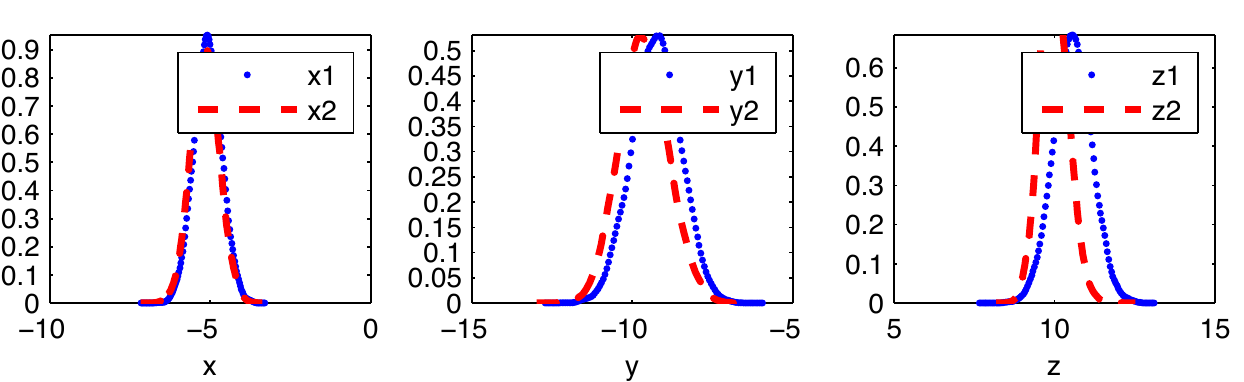}
 \caption{Linear measurement: Comparison posterior for LBU ($n=1$) and 
 NLBU ($n=2$) after second update}
\label{F:exp-NB-6}
\end{figure}

As a first set of experiments we take the measurement operator to be
linear in $q$; $Y(q) = q =(x,y,z)$, i.e.\ we can observe the \emph{whole}
state directly.  At the moment we consider updates after each day---whereas
in \refig{exp-B-2} the updates were performed every 10 days.
The results for the pdfs of the state variables are shown in \refig{exp-NB-5},
where the prior and the posterior pdf for a LBU after one update
are given.  Then we do the same experiment, but with a \emph{quadratic} nonlinear
BU (NLBU) with $n=2$.  The results for the posterior pdfs
are given in \refig{exp-NB-5.5}, where the
linear update is dotted in blue, and the full red line is the quadratic NLBU; there is
hardly any difference between the two.  This might have been expected after
our discussion at the end of \refSS{prior}.  If we go on to the second update---after
two days---some differences appear, the results for the posterior pdfs
are in \refig{exp-NB-6}.

As the differences between LBU and NLBU with $n=2$ were small---we take this as an
indication that the LBU is not too inaccurate an approximation to the conditional
expectation---we change the experiment and
take a nonlinear measurement function, which is now cubic: $Y(q) =(x^3,y^3,z^3)$.
As discussed at the end of \refSS{prior}, we now expect larger differences
between LBU and NLBU.
%
\begin{figure}[!ht]
\centering
 \includegraphics[width=0.9\textwidth,height=0.2\textheight]{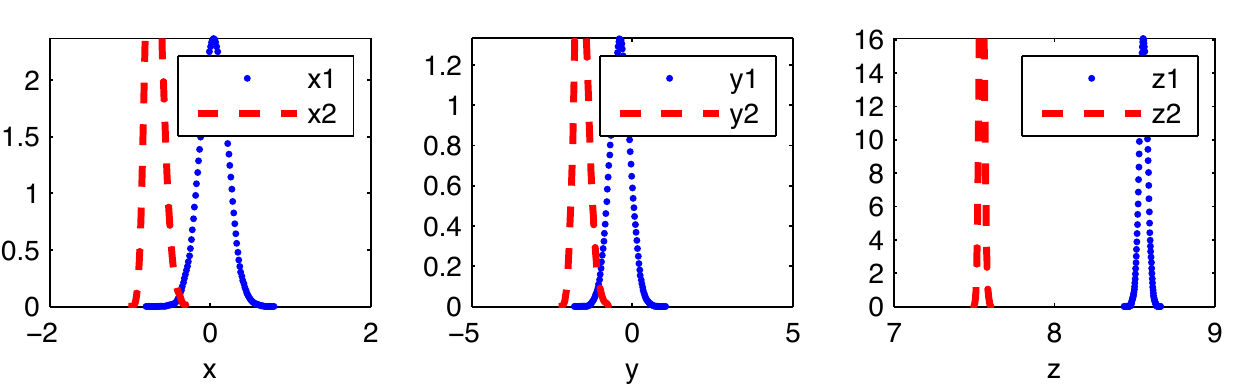}
 \caption{Cubic measurement: Comparison posterior for LBU ($n=1$) and 
 NLBU ($n=2$) after one update}
\label{F:exp-NB-9}
\end{figure}

These differences in posterior pdfs after one update may be gleaned 
from \refig{exp-NB-9}, and they are indeed larger
than in the linear case \refig{exp-NB-5.5}, due to the strongly nonlinear
measurement operator.
\begin{figure}[!ht]
\centering
 \includegraphics[width=0.9\textwidth,height=0.3\textheight]{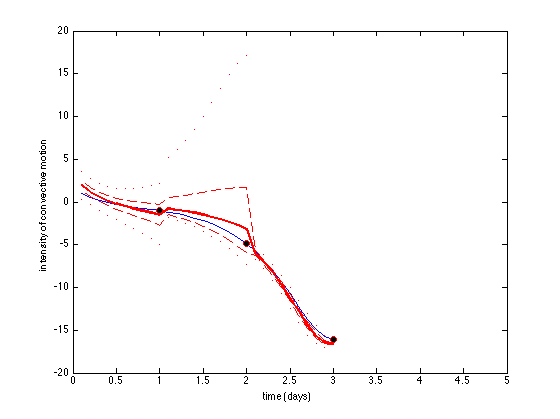}
 \caption{Partial state trajectory with uncertainty and three updates}
\label{F:exp-NB-9.5}
\end{figure}

As the cubic is quite a strong nonlinearity, we performed a set of experiments
where the measurement function is $Y(q) = (x |x|,y |y|,z |z|)$;
only a quadratic nonlinearity, but no loss of information about the sign like
in the small example at the end of \refSS{prior}.  The updates are performed every
day, the \refig{exp-NB-9.5}, which shows the trajectory of one state variable,
corresponds in that way to \refig{exp-B-2}.
\begin{figure}[!ht]
\centering
 \includegraphics[width=0.9\textwidth,height=0.2\textheight]{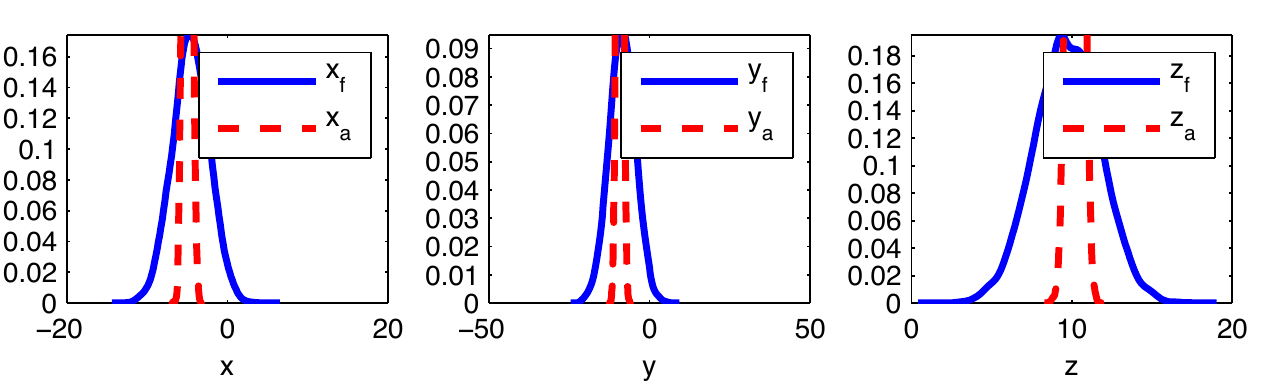}
 \caption{Quadratic measurement: Comparison posterior for LBU ($n=1$) and 
 NLBU ($n=2$) after one update}
\label{F:exp-NB-10}
\end{figure}

The results for the $2$-nd update are displayed for the posterior pdfs
in \refig{exp-NB-10}.  This has to be compared \refig{exp-NB-6}, and the
differences are indeed much larger.

%
%
%
%
%
%
%


%

\section{Conclusion} \label{S:concl}
Here we have tried to show the connection between inverse problems and
uncertainty quantification.  An abstract model of a system was introduced, together
with a measurement operator, which provides a possibility to predict---in a
probabilistic sense---a measurement.  The framework chosen is that of Bayesian
analysis, where uncertain quantities are modelled as random variables.
New information leads to an update of the probabilistic description via
Bayes's rule.

After elaborating on the---often not well-known---connection between conditional
probabilities as in Bayes's rule and conditional expectation, we set out to
compute and---necessarily---approximate the conditional expectation.  As a
polynomial approximation as chosen, there is the choice up to which degree
one should go.  The case with up to linear terms---the linear Bayesian update---is
best known and intimately connected with the well-known Kalman filter.  In addition,
we show how to compute approximations of higher order.

There are several possibilities on how one may choose a numerical realisation of
these theoretical concepts, and we decided on functional or spectral approximations.
It turns out that this approach goes very well with recent very efficient approximation
methods building on separated or so-called low-rank tensor approximations.

Starting with the linear Bayesian update, we show a series of examples of increasing
complexity.  The method works well in all cases.  One of the examples is then chosen
to show the nonlinear Bayesian update, where we go up to quadratic terms.  A series
of experiments is chosen with different measurement operators, which have quite
a marked influence on whether the linear and quadratic update are close to
each other.

%
%
%
%
%




\bibliography{\thebib/jabbrevlong,\thebib/matthies_BU_paper-1,\thebib/phys_D,\thebib/fa,\thebib/risk,\thebib/fuq-new,\thebib/highdim}


{ 
   \tiny
       \texttt{\RCSId} 
   }



\end{document}